\documentclass[11pt,letterpaper]{amsart}

\usepackage{amsmath,amssymb,amsthm,verbatim, csquotes,mathrsfs,stmaryrd,amsfonts}
\usepackage{hyperref}
\usepackage{xcolor}
\usepackage{esint}
\usepackage{mathtools}
\usepackage{graphicx}
\usepackage{float}
\usepackage{caption}
\mathtoolsset{showonlyrefs}
\usepackage{geometry}

\usepackage{bm}
\usepackage{pdfpages}
\usepackage{import}
\usepackage[shortlabels]{enumitem}

\newtheorem{theorem}{Theorem}[section]
\newtheorem{lemma}[theorem]{Lemma}
\newtheorem{proposition}[theorem]{Proposition}
\newtheorem{definition}[theorem]{Definition}
\newtheorem{corollary}[theorem]{Corollary}
\newtheorem{remark}[theorem]{Remark}

\numberwithin{equation} {section}

\newcommand{\R}{{\mathbb R}}
\renewcommand{\O}{{\mathbb O}}
\newcommand{\SO}{{\mathbb {SO}}}

\begin{document}

\title[]{Curvature varifolds with given oriented boundary}
\date{\today}

\author[Kuwert]{Ernst Kuwert}
\address[E.K]{Mathematisches Institut\\
Universit\"at Freiburg\\
Ernst-Zermelo-Str.1\\
79104\\
\newline Freiburg\\ Germany}
\email{ernst.kuwert@math.uni-freiburg.de}

\author[Zhang]{Xuwen Zhang}
	\address[X.Z]{Mathematisches Institut\\
		Universit\"at Freiburg\\
	Ernst-Zermelo-Str.1\\
		79104\\
  \newline Freiburg\\ Germany}
\email{xuwen.zhang@math.uni-freiburg.de}

\begin{abstract}
Let $\gamma\subset\mathbb R^n$ be a closed, embedded, connected, oriented $C^1$ submanifold in $\mathbb R^n$ of dimension $m-1$ where $2\leq m<n$.
We prove that the area of any immersed $m$-dimensional submanifold with oriented boundary $\gamma$ is bounded in terms of its curvature energy, unless $\gamma$ lies in an $m$-plane.
The problems leads to a weak formulation involving oriented curvature varifolds with given oriented boundary, and we solve a boundary-value minimization problem with obstacle in this class.
\\

\end{abstract}

\maketitle
\section{Introduction}

Let $\gamma\subset\mathbb R^n$ be a closed, embedded, connected, oriented $C^1$ submanifold in $\mathbb R^n$ of dimension $m-1$ where $2\leq m<n$.
We denote by $\mathbf S^m(\gamma)$ the class of compact, oriented, $C^1$ immersions $\Sigma\hookrightarrow\mathbb R^n$ of dimension $m$, such that the restriction to $\partial \Sigma$ parametrizes $\gamma$, or more generally, is an oriented covering of $\gamma$. Assuming that the immersion is of class $C^2$ in the interior we can then consider the $L^p$ curvature energy
$$
\int_\Sigma |A|^p\,d\mathcal H^m
\quad \mbox{ where }1 \leq p < \infty.
$$
Here $\Sigma$ is endowed with the induced metric, and $A$ is the second fundamental form.

In general, the area of an immersed submanifold with a given oriented boundary $\gamma$ is not bounded in terms of its curvature energy. For example, a $2$-disk $D_1:=\{|x|^2=1\}$ in $\mathbb R^3$ has positive area but zero curvature.
Concerning an immersed example, one can simply take an immersed annulus into a plane in $\mathbb{R}^3$, given by
$$
F:[a,b]\times\mathbb S^1\to\mathbb R^3,\quad
F(r,\theta)=(r\cos{k\theta},r\sin{k\theta},0),
$$
where $0<a<b$ and $k\in\mathbb N$.
As $k\to\infty$ the immersed annulus has area goes to infinity, but curvature always zero.
These examples show that it is necessary for $\gamma$ to be ``non-planar'' such that the area could be bounded in terms of the curvature energy.

As a first main result, we prove such a bound for submanifolds with given oriented boundary in any dimension and codimension.


\begin{theorem}\label{Thm:area-bounds}
Let $\Sigma\in\mathbf{S}^m(\gamma)$ satisfy ${\rm diam}\Sigma\leq R<\infty$.
If $\gamma$ is not contained in any affine $m$-plane,
then for every $p\in(1,\infty)$, there exists $C=C(m,n,p,\gamma,R)>0$ such that
\begin{align}\label{ineq:area-bounds}
|\Sigma|+|\partial\Sigma|
\leq C\int_\Sigma|A|^p\,d\mathcal H^m.
\end{align}
\end{theorem}
Let us make some comments.
First, note that both the assumptions that $\gamma$ is not contained in any affine $m$-plane, and that diameter of $\Sigma$ is bounded are necessary. For example, consider the spherical caps $S^+_r:=\left\{x=(x_1,x_2,x_3)\in\mathbb R^3:|x-(0,0,\sqrt{r^2-1})|=r, x_3\geq0\right\}$ ($r>1$) with the same boundary curve given by the unit circle $\gamma_1:=\{x\in\mathbb R^3:|x|=1, x_3=0\}$. As $r\to\infty$ we have $|S^+_r|\to\infty$, but $\int_{S^+_r}|A|^2\,d\mathcal H^2\to8\pi$.
By a simple perturbation of the family of caps $\{S^+_r\}_{r>1}$, we see that the area bound fails to be true in this case.
Second, it is interesting to note that, since $\Sigma$ is an isometric immersion and $\partial\Sigma$ is an oriented covering of $\gamma$, the bounds \eqref{ineq:area-bounds} in fact means that we can estimate the degree of the covering directly by the curvature energy.
Third, in the case $m=2<n$, it is well-known to experts that as a direct application of the boundary version of Simon's monotonicity formula \cite{Simon93,KS04,NP20}, below the (critical) energy gap $\int_\Sigma|A|^2\,d\mathcal H^2<8\pi$, immersions $\Sigma\in\mathbf{S}^2(\gamma)$ have a priori a diameter bound.

There are a number of papers dealing with the relation between area and curvature. Here we list those closely related to our work.
By embedding a Riemannian manifold into an ambient Euclidean space then using Hutchinson's curvature varifolds \cite{Hutchinson86}, Mondino \cite{Mondino14} proved by a contradiction argument the mass bounds for $m$-varifolds in a Riemannian manifold in terms of their curvature, under the condition that there exist no curvature varifolds with zero curvature on the given Riemannian manifold.
In the special case $m=2$, Bangert-Kuwert \cite{BK25} showed that for such a mass bound to hold, it is sufficient to rule out smooth totally geodesic immersions.

For surfaces with boundary, Novaga-Pozzetta \cite{NP20} studied the minimization of Willmore energy in the context of Mantegazza's curvature varifolds with boundary \cite{Mantegazza96}.
As an intermediate step, they proved that if the Willmore energy is below the energy threshold $\int|H|^2\,d\mathcal H^2<4\pi$, then $|\Sigma|$ is controlled in terms of the Willmore energy $\int|H|^2\,d\mathcal H^2$ and the boundary length $|\gamma|$, see also \cite{Pozzetta21} for such an estimate.
More recently, Kuwert-M\"uller \cite{KM22} proved that for any given container $\Omega$, if the intersection of any affine $m$-plane with $\Omega$ is not orthogonal, then any immersion into $\Omega$ with orthogonal boundary along $\partial\Omega$ satisfies a bound of the form \eqref{ineq:area-bounds}.

In our case, we prove Theorem \ref{Thm:area-bounds} using a contradiction argument aligned with the works \cite{Mondino14,KM22}. A key point is that one should introduce a suitable weak formulation of immersions with given oriented boundary such that a nice compactness theorem can be established and used for the contradiction argument.
In our case, we work with the so-called oriented varifolds, which are Radon measures on the oriented Grassmannian bundle $\mathbb R^n\times\widetilde G(m,n)$.
To our goal, we introduce the following subspace of $\gamma\times\widetilde G(m,n)$, which provides all possible choices of oriented tangent planes at the boundary points for any $\Sigma\in\mathbf{S}^m(\gamma)$:
\begin{align}\label{defn:G^o(m,n,gamma)(x)}
\widetilde G_{m,n,\gamma}
:=\left\{(p,\omega)\in\gamma\times\widetilde G(m,n):
v\wedge\omega=0\text{ for every }v\in T_p\gamma\right\}.
\end{align}
Such a definition is in spirit aligned with the works \cite{KM22,WZ25,Zhang26}, and as in these works, we can define the {\em weak unit inwards co-normal} as
\begin{align}\label{defn:eta(x,om)}
\eta(p,\omega)
:=(-1)^{n-m}\star\left(\tau_\gamma(p)\wedge\star\omega\right),\quad\forall(p,\omega)\in\widetilde G_{m,n,\gamma},
\end{align}
where $\star$ is the Hodge star operator, and $\tau_\gamma(p)$ is the chosen oriented unit simple $(m-1)$-vector associated to the oriented tangent plane of $\gamma$ at $p$.
One can readily check by definition that $\eta(p,\omega)$ is a unit vector in $\mathbb R^n$, tangent to the $m$-plane spanned by $\omega$ and perpendicular to $T_p\gamma$.
Moreover, for $V=\mathcal H^m\llcorner\Sigma\otimes\delta_{T_p\Sigma}$ and $\Gamma=\mathcal H^{m-1}\llcorner\gamma\otimes\delta_{T_p\gamma}$ induced by $\Sigma\in{\mathbf S}^m(\gamma)$, the first variation formula (cf. \cite[Lemma 39.4]{Simon83}) holds: for any $X\in C^1_c(\mathbb R^n;\mathbb R^n)$,
\begin{align}\label{eq:1st-variation}
\delta V[X]
=\int_\Sigma{\rm div}X\,d\mathcal H^m
=&-\int\left<H,X\right>d\mathcal H^m-\int\left<\eta(p,\omega),X(p)\right>d\Gamma(p,\omega)\\
=&-\int\left<H,X\right>d\mathcal H^m-\int_\gamma\left<\eta(p,T_p\gamma),X(p)\right>d\mathcal H^{m-1}(p),
\end{align}
where $H$ is the mean curvature vector and $\eta(p,T_p\gamma)$ agrees with the classical inwards unit co-normal along $\gamma=\partial\Sigma\subset\Sigma$.
Now we can define the class of oriented curvature varifolds with given oriented boundary $\gamma$, say $\widetilde{\bf CV}_{m,\gamma}(\mathbb R^n)$, by allowing the test functions $X$ in \eqref{eq:1st-variation} depend also on the Grassmannian variable, as done in \cite{Hutchinson86}.

We will show in Proposition \ref{Prop:initial-cptness} that, by compactness of $\widetilde{\bf CV}_{m,\gamma}(\mathbb R^n)$, the given oriented boundary condition carries over to the limits. In Theorem \ref{Thm:structure-V-Gamma-S=0}, we prove that there is no non-zero curvature varifold in $\widetilde{\bf CV}_{m,\gamma}(\mathbb R^n)$ with vanishing curvature and bounded support, provided that $\gamma$ is not contained in any $m$-plane.
Combining these two facts, we obtain Theorem \ref{Thm:area-bounds} even in the weak class $\widetilde{\bf CV}_{m,\gamma}(\mathbb R^n)$.
Moreover, as a consequence, we can show the existence of minimizers of the following boundary value problem with obstacle among the integer rectifiable subclass of $\widetilde{\bf CV}_{m,\gamma}(\mathbb R^n)$, say $\widetilde{\bf ICV}_{m,\gamma}(\mathbb R^n)$:
\begin{theorem}\label{thm:minimizer-exis}
Given $p\in(1,\infty)$, $R\in(0,\infty)$ and a positive integer $k$, the following infimum is attained whenever the class is non-empty:
\begin{align}\label{defn:mmmzer-L^p}
&\beta^p_{k,R}(\gamma)\\
:=&\inf\left\{\int|\widetilde B|^p\,dV:V\in\widetilde{\bf ICV}_{m,\gamma}(\mathbb R^n), {\rm spt}\|V\|\subset\overline{B_R(0)}, \pi_\#\Gamma=k\mathcal H^{m-1}\llcorner\gamma\right\}.
\end{align}
Moreover, if $\gamma$ is not contained in any affine $m$-plane, then $\beta^p_{k,R}(\gamma)>0$.
\end{theorem}
Here $V$ denotes the weak immersion, $\widetilde B$ is its weak curvature and $\Gamma$ is its weak boundary. $\pi_\#\Gamma=k\mathcal H^{m-1}\llcorner\gamma$ means that the weak immersion has boundary covering $\gamma$ $k$-times.
We note that when minimizing \eqref{defn:mmmzer-L^p}, the co-normals of the weak surfaces are allowed to change freely on $\gamma$, and if the minimizer of \eqref{defn:mmmzer-L^p} were smooth, then the Euler-Lagrange equation is in fact a fourth order equation in the interior, see e.g. \cite{KLL15}, and a second order equation on the boundary. In the minimization problem of Willmore energy, the same boundary condition has been considered by Pozzetta in \cite{Pozzetta21}, and a Dirichlet-type boundary was studied in \cite{NP20}.

It is worth mentioning that the minimization problem \eqref{defn:mmmzer-L^p} cannot be solved using Mantegazza's unoriented curvature varifold with boundary \cite{Mantegazza96}. More precisely, in his definition, the boundary of an integral curvature varifold $V$ is given by a vector-valued Radon measure, say $\partial V$, and the boundary condition should be formulated as $|\partial V|=k\mathcal H^{m-1}\llcorner\gamma$, where $|\partial V|$ is the total variation of $\partial V$.
However, we note that there is a lack of nice compactness in such settings. For example, consider the union $C_r:=S_r^+\cup S_r^-$ ($r>1$), where $S_r^+$ is the spherical cap defined above, and $S_r^-$ is the reflection of $S_r^+$ across the hyperplane $\{x_3=0\}$.
One can show that $C_r$ induces an integral curvature varifold $V_r$ with boundary $\partial V_r$, such that $|\partial V_r|=2\mathcal H^{m-1}\llcorner\gamma_1$ (recall that $\gamma_1$ is the unit circle).
However, letting $r\searrow1$, $(V_r,\partial V_r)$ converges as curvature varifold with boundary to $(V_1,\partial V_1)$, where $V_1$ is induced by the unit sphere $\{x\in\mathbb R^3:|x|^2=1\}$ and $\partial V_1=0$, namely, the boundary condition fails to be carried over to the limit.

\

\noindent{\em The paper is organized as follows.}
In Section \ref{sec:2} we collect some basic facts from Geometric measure theory and study the structure of oriented varifolds.
In Section \ref{sec:3} we define the weak class of submanifolds that we work with, namely, the oriented curvature varifolds with given oriented boundary, and study its basic properties. Then we consider the zero curvature case in Section \ref{sec:4} and prove the key structural result (Theorem \ref{Thm:structure-V-Gamma-S=0}).
Section \ref{sec:5} is devoted to establishing the compactness result (Proposition \ref{Prop:initial-cptness}) for oriented curvature varifolds with given oriented boundary, which is essentially based on the key estimate Lemma \ref{Lem:bdry-mass-estimate}.
In Section \ref{sec:6} we prove the area bound (Theorem \ref{Thm:area-bounds}) and study in detail the two dimensional case in Section \ref{sec:6-mass-topo-bds}.
In Section \ref{sec:7} we study the minimization problem \eqref{defn:mmmzer-L^p} and prove Theorem \ref{thm:minimizer-exis}.

\

\noindent{\em Acknowledgments.} The authors thank Fabian Rupp for helpful discussions on this project.
This work is supported by the DFG through the Grant IRTG3132 “Analysis of Nonlinear PDE”.

\

\noindent{\em AI disclosure.} AI was not used in the development of this paper.

\section{Preliminaries}\label{sec:2}

We first collect two useful facts from Geometric measure theory.
The first one is the well-known Constancy Theorem.
\begin{lemma}\label{Lem:Constancy}
Let $\mu$ be a Radon measure on an open domain $U\subset\mathbb{R}^m$, and assume that for all $f\in C^\infty_c(U)$
\begin{align}
\int_U\partial_\alpha f(x)\,d\mu(x)=0,\quad\forall\alpha=1,\cdots,m.
\end{align}
Then $\mu=\theta\mathcal{L}^m$ on $U$ for some $\theta\in[0,\infty)$.
\end{lemma}
\begin{proof}

See \cite[Theorem 26.27]{Simon83}.
\end{proof}

The second one is disintegration for Radon measures, see e.g. \cite{AGS08,Simon83}.
\begin{lemma}\label{lem:disintegration}
Let $\mu$ be a Radon measure on $X\times Y$ with compact support, where $X,Y$ are metric spaces, and denote by $\pi:X\times Y\rightarrow X$ the projection map.
Then $\|\mu\|=\pi_\#\mu$ is a Radon measure, and there is a $\|\mu\|$-a.e. uniquely determined family $(\mu^x)_{x\in X}$ of Radon probability measures on $Y$, such that for any Borel function $\phi:X\times Y\to[0,\infty]$ one has
\begin{align}
\int_{X\times Y}\phi \,d\mu
=\int_X\int_Y\phi(x,y)\,d\mu^x(y){d}\|\mu\|(x).
\end{align}
\end{lemma}

\subsection{Oriented \texorpdfstring{$m$}{m}-varifolds in \texorpdfstring{$\mathbb{R}^{n}$}{R{n}}}


\subsubsection{Grassmannian}
The unoriented and oriented Grassmannians are 
$$
G(m,n)={\rm V}_m(\mathbb{R}^n)/\O(m)
\quad \mbox{ and } \quad
\widetilde{G}(m,n)={\rm V}_m(\mathbb{R}^n)/\SO(m), 
$$
where ${\rm V}_m(\mathbb{R}^n) =
\{X \in \mathbb{R}^{n\times m}:X^{T}X = I_m\}$ 
is the Stiefel manifold of orthonormal $m$-frames 
in $\R^n$. We have the obvious covering 
\begin{align*}
q:\widetilde{G}(m,n)\to G(m,n),
\quad
[X]_{SO(m)}\mapsto[X]_{O(m)}.
\end{align*}
On $V_m(\R^n)$ we have the metric induced from $\R^{n \times m}$, 
and the metrics on the Gra{\ss}mannians are chosen such that
the quotient projections are Riemannian submersions.
For either $[X]_{O(m)}\in G(m,n)$ or $[X]_{SO(m)}\in \widetilde{G}(m,n)$,
the tangent space at $[X]$ is identified with the horizontal space 
${\mathcal H}_X = \{Z \in \R^{n \times m}: X^T Z=0\}$.
The quotient metric is given by 
\begin{align*}
g_{[X]}(Z,W)={\rm tr}(Z^{T}W) \quad \mbox{ where }
Z,W \in {\mathcal H}_X.
\end{align*}
Now consider the map
\begin{align*}
\iota:G(m,n)\to {\rm Sym}_n(\mathbb{R}),\quad
[X]_{O(m)}\mapsto XX^{T}.
\end{align*}
This is well defined, since for every $A\in O(m)$ one has
\begin{align*}
(XA)(XA)^{T}=XAA^{T}X^{T}=XX^{T}.
\end{align*}
Its image is the space of rank $m$ orthogonal projection matrices:
\begin{align*}
\iota(G(m,n))
=\mathcal{P}_{m,n}
:=\left\{
P\in{\rm Sym}_n(\mathbb{R}):
P^{2}=P,\ {\rm tr}(P)=m
\right\}.
\end{align*}
We equip ${\rm Sym}_n(\mathbb{R})$ with the inner product
$\left< A,B\right>_{{\rm Sym}}
=
\frac{1}{2}{\rm tr}(AB)$, and note that the differential of $\iota$ at $[X]_{O(m)}$ is
\begin{align*}
d\iota_{[X]_{O(m)}}(Z)=ZX^{T}+XZ^{T}.
\end{align*}
Using \(X^{T}Z=X^{T}W=0\), we thus obtain
\begin{align*}
\left< d\iota_{[X]}(Z),d\iota_{[X]}(W) \right>_{{\rm Sym}}
&=\frac{1}{2}{\rm tr}\left((ZX^{T}+XZ^{T})(WX^{T}+XW^{T})\right)\\
&={\rm tr}(Z^{T}W)
=g_{[X]_{O(m)}}(Z,W),
\end{align*}
which shows that $\iota:G(m,n)\to {\rm Sym}_n(\mathbb{R})$
is an isometric embedding.

For the oriented Grassmannian, define
\begin{align*}
\widetilde{\iota}:
\widetilde{G}(m,n)
\to
\Lambda_{m}(\mathbb{R}^{n}),\quad
[v_1,\ldots,v_m]_{SO(m)}\mapsto
v_1\wedge\cdots\wedge v_m,
\end{align*}
where $(v_1,\ldots,v_m)$ is an oriented orthonormal frame. $\widetilde\iota$ is
well defined, since a change of $(v_1,\ldots,v_m)$ by some $A\in SO(m)$
multiplies the wedge product by $\det A=1$.
Its image is the space of unit simple $m$-vectors:
\begin{align*}
\widetilde\iota(\widetilde G(m,n))
=\widetilde{\mathcal{G}}_{m,n}
:=\left\{v_1\wedge\cdots\wedge v_m\in\Lambda_m(\mathbb{R}^n):(v_1,\cdots,v_m)\text{ is an orthonormal }m\text{-frame}\right\}.
\end{align*}
We equip $\Lambda_{m}(\mathbb{R}^{n})$ with the Euclidean inner product
induced from that of \(\mathbb{R}^{n}\), and note that for $X=(v_1,\ldots,v_m)$ and
$Z=(z_1,\ldots,z_m)$ with $z_i\in{\rm span}(v_1,\ldots,v_m)^{\perp}$, there holds
\begin{align*}
d\widetilde{\iota}_{[X]_{SO(m)}}(Z)
=
\sum_{i=1}^{m}
v_1\wedge\cdots\wedge z_i\wedge\cdots\wedge v_m,
\end{align*}
we thus deduce
\begin{align*}
\left<
d\widetilde{\iota}_{[X]_{SO(m)}}(Z),
d\widetilde{\iota}_{[X]_{SO(m)}}(W)
 \right>_{\Lambda_{m}(\mathbb{R}^n)}
=\sum_{i=1}^{m}\left< z_i,w_i\right>
={\rm tr}(Z^{T}W)
=g_{[X]_{SO(m)}}(Z,W).
\end{align*}
Hence $\widetilde{\iota}:
\widetilde{G}(m,n)\to\Lambda_{m}(\mathbb{R}^{n})$
is also an isometric embedding.

Finally, the covering $q:\widetilde G(m,n)\to G(m,n)$ induces a covering between the two
embeddings:
We define
\begin{align*}
\mathbf{q}:\widetilde{\mathcal{G}}_{m,n}
\to
\mathcal{P}_{m,n},\quad
v_1\wedge\cdots\wedge v_m
\mapsto
\sum_{i=1}^{m}v_iv_i^{T},
\end{align*}
where \((v_1,\ldots,v_m)\) is an orthonormal frame.
Clearly, the maps satisfy
$\mathbf{q}\circ\widetilde{\iota}=\iota\circ q$.

The involution map on $\widetilde{\mathcal{G}}_{m,n}$ is given by
\begin{align*}
\mathcal{I}:\widetilde{\mathcal G}_{m,n}
\to
\widetilde{\mathcal G}_{m,n},
\quad
v_1\wedge\cdots\wedge v_m
\longmapsto
-v_1\wedge\cdots\wedge v_m,
\end{align*}
which satisfies
$\mathcal{I}^2={\rm Id}_{\widetilde{\mathcal G}_{m,n}}$, and it is easy to check that $\mathcal{I}$ is an isometry.
Note also $\mathbf q\circ\mathcal I=\mathbf q$, while for every $\omega\in\widetilde{\mathcal G}_{m,n}$, $\mathbf q^{-1}\left(\mathbf q(\omega)\right)=\{\omega,\mathcal I(\omega)\}=\{\omega,-\omega\}$.
Hence, $\mathcal P_{m,n}
\cong
\widetilde{\mathcal G}_{m,n}/\left<\mathcal I\right>$,
where the quotient is equipped with the quotient Riemannian metric, and $\mathbf{q}:\widetilde{\mathcal{G}}_{m,n}
\to
\mathcal{P}_{m,n}$ is a Riemannian double covering.

\begin{remark}
\normalfont
For the rest of the paper, we identify $G(m,n),\widetilde G(m,n)$ with their isometric embeddings $\mathcal{P}_{m,n}$, $\widetilde{\mathcal{G}}_{m,n}$ into Euclidean spaces ${\rm Sym}_n(\mathbb{R})$, $\Lambda_m(\mathbb{R}^n)$.
\end{remark}

\subsubsection{Oriented varifolds and its structure}\label{Sec:oriented-structure}

For any open $U\subset\mathbb{R}^n$, an \textit{oriented $m$-varifold on $U$} is a Radon measure $V$ on $\widetilde{G}_{m,n}(U):=U\times \widetilde{G}(m,n)$.
With a slight abuse of notation, we also denote by $\mathbf{q}$ the map from $\widetilde{G}_{m,n}(U)$ to $G_{m,n}(U):=U\times G(m,n)$, sending $(x,\omega)$ to $(x,\mathbf{q}(\omega))$; and by $\mathcal{I}$ the isometric involution from $\widetilde{G}_{m,n}(U)$ to $\widetilde{G}_{m,n}(U)$, sending $(x,\omega)$ to $(x,\mathcal{I}(\omega))$.

We write $V=V^s+V^{as}$ to be its \textit{symmetric} and \textit{anti-symmetric} part respectively, where
\begin{align}\label{defn:V^s-V^as}
V^s
:=\frac12(V+\mathcal{I}_\#V),\quad\text{ and }
V^{as}
:=\frac12(V-\mathcal{I}_\#V).
\end{align}
Note that $V^s$ is a Radon measure, while $V^{as}$ is only a signed Radon measure.

We put further $V^0=\mathbf{q}_\#V$. Since $\mathbf{q}\circ\mathcal{I}=\mathbf{q}$, thus we infer
\begin{align}\label{defn:V^0}
\mathbf{q}_\#V^s
=V^0,\quad\text{ and }
\mathbf{q}_\#V^{as}
=0.
\end{align}

By \eqref{defn:V^s-V^as}, it is easy to see that $V$ is absolutely continuous with respect to $V^s$, and hence $V^{as}$ is also absolutely continuous with respect to $V^s$ due to $V^{as}=V-V^{s}$.
By Radon-Nikodym theorem, there exists a $V^s$-measurable function $\sigma:\widetilde{G}_{m,n}(U)\rightarrow\mathbb{R}$, such that
\begin{align}\label{defn:V^as-sigma}
V^{as}(\phi)
=\int\phi(x,\omega)\sigma(x,\omega)\,{d} V^s(x,\omega),\quad\forall\phi\in C^0_c(U\times \Lambda_m(\mathbb{R}^n)).
\end{align}

\begin{proposition}\label{Prop:sigma-odd}
Under the above notations, $\sigma$ is an odd function on the second variable and takes values in $[-1,1]$ for $V^s$-a.e. $(x,\omega)$.
\end{proposition}
\begin{proof}
First by \eqref{defn:V^s-V^as}, we note that
\begin{align*}
\mathcal{I}_\#V^s=V^s,\quad
\mathcal{I}_\#V^{as}=-V^{as}.
\end{align*}
Hence by definition of push-forward, we find $V^{as}(\phi\circ\mathcal{I})=(\mathcal{I}_\#V^{as})(\phi)=-V^{as}(\phi)$.

On the other hand, by definition of push-forward of measures we have
\begin{align*}
V^{as}(\phi\circ\mathcal{I})
\overset{\eqref{defn:V^as-sigma}}{=}&\int(\phi\circ\mathcal{I})\sigma\,{d}V^s
=\int\left(\phi(\sigma\circ\mathcal{I})\right)\circ\mathcal{I}\,{d}V^s\\
=&\int\phi\left(\sigma\circ\mathcal{I}\right)\,{d}(\mathcal{I}_\#V^s)
=\int\phi(\sigma\circ\mathcal{I})\,{d} V^s.
\end{align*}
Combining we obtain for all $\phi\in C^0_c(U\times\Lambda_m(\mathbb{R}^n))$
\begin{align}
\int\phi\left(\sigma+\sigma\circ\mathcal{I}\right)\,{d} V^s
=0.
\end{align}
By a density argument we conclude that $\sigma\circ\mathcal{I}=-\sigma$, i.e., $\sigma$ is odd on the second variable for $V^s$-a.e. $(x,\omega)$.

To show $\sigma(x,\omega)\in[-1,1]$ for $V^s$-a.e. $(x,\omega)$, it suffices to observe that
\begin{align*}
V=V^s+V^{as}
=(1+\sigma)V^s,\quad
\mathcal{I}_\#V
=V^s-V^{as}
=(1-\sigma)V^s
\end{align*}
are both non-negative measures, and hence
\begin{align*}
-1\leq\sigma\leq1\text{ for }V^s\text{-a.e.}
\end{align*}
\end{proof}

\section{Oriented curvature varifolds with given oriented boundary}\label{sec:3}

\subsection{Boundary condition}

Now we generalize smooth surfaces with given oriented boundary to the context of oriented varifolds.

Let $\gamma\subset\mathbb R^n$ be a closed connected oriented embedded $C^1$ submanifold in $\mathbb R^n$ of dimension $m-1$. We use the induced Euclidean metric and regard $T_p\gamma$ as a subspace of $\mathbb R^n$ and denote by $\tau_\gamma(p)\in\Lambda_{m-1}(\mathbb R^n)$ the oriented unit simple $(m-1)$-vector associated to $T_p\gamma$.

The oriented Grassmannian bundle associated with the given oriented boundary $\gamma$ is defined by $\widetilde G_{m,n,\gamma}$ as in \eqref{defn:G^o(m,n,gamma)(x)}.
We endow $\widetilde G_{m,n,\gamma}$ with the subspace topology inherited
from $\gamma\times\widetilde G(m,n)$ and denote its fiber as $\widetilde G_{m,n,\gamma}(x):=\left\{\omega\in\widetilde G(m,n):v\wedge\omega=0\text{ for every }v\in T_x\gamma\right\}$.
We also use the interior product $\iota_v\omega$ given by
\begin{align*}
\left<\iota_v\omega,\alpha\right>
:=\left<\omega,v\wedge\alpha\right>,
\quad\forall \alpha\in\Lambda_{m-1}(\mathbb R^n).
\end{align*}
For every $(p,\omega)\in\widetilde G_{m,n,\gamma}$, define the {\em weak inner unit conormal} $\eta(p,\omega)$ as in \eqref{defn:eta(x,om)}.
By the definition of $\widetilde G_{m,n,\gamma}$, the subspaces represented by $\tau_\gamma(p)$ and $\star\omega$ are orthogonal. Hence $\tau_\gamma(p)\wedge\star\omega$
is a unit simple $(n-1)$-vector, so that $\eta(p,\omega)$ is a unit vector.
More precisely, $\eta(p,\omega)\in{\rm span}(\omega)\cap(T_p\gamma)^\perp$, and it satisfies 
\begin{align*}
\tau_\gamma(p)\wedge\eta(p,\omega)=\omega,
\end{align*}
thanks to \eqref{defn:eta(x,om)}.
We also note that it can be characterized by
\begin{align}\label{eq:immersed-eta-contraction}
\iota_{\eta(p,\omega)}\omega
=(-1)^{m-1}\tau_\gamma(p).
\end{align}

\subsection{Oriented curvature varifolds with given oriented boundary}
In the spirit of \cite{Hutchinson86} (see also \cite[Section 2]{KM22}), we first derive the curvature identity for smooth surfaces.

Let $\gamma\subset U\subset\mathbb R^n$ be as before, let $\Sigma$ be a compact, oriented manifold of dimension $m$ with nonempty boundary, and let $F:\Sigma\to U$ be a $C^1$ immersion that is $C^2$ in the interior. We equip $\Sigma$ with the metric induced by $F$ and identify every tangent space of $\Sigma$ and $\partial\Sigma$ with its image in $\mathbb R^n$ under the pushforward of $F$, taking derivatives along these spaces are also understood through this identification.

Now we denote by $\omega_\Sigma(z)\in\Lambda_m(\mathbb R^n)$ the oriented unit simple $m$-vector associated to $T_z\Sigma$, and $\eta_{\partial\Sigma}(z)\in\mathbb R^n$ the inward-pointing unit conormal.
The boundary orientation is given by
\begin{align}\label{eq:tau_partial-Sigma-and-contraction}
\tau_{\partial\Sigma}(z)
:=(-1)^{m-1}\iota_{\eta_{\partial\Sigma}(z)}\omega_\Sigma(z).
\end{align}
For any $\Sigma\in\mathbf{S}^m(\gamma)$, since $F\mid_{\partial\Sigma}$ is an orientation-preserving covering of $\gamma$,  we have
\begin{align}\label{eq:boundary-orientation-compatibility}
\tau_{\partial\Sigma}(z)=\tau_\gamma(F(z))\quad\text{for every }z\in\partial\Sigma.
\end{align}
Moreover, since $\partial\Sigma$ is compact and $\gamma$ is connected, we know this covering is finite and $F(\partial\Sigma)=\gamma$.
For any $z\in\partial\Sigma$, we have $T_{F(z)}\gamma=T_z\partial\Sigma\subset T_z\Sigma$, thus by definition $(F(z),\omega_\Sigma(z))\in\widetilde G_{m,n,\gamma}$. By \eqref{eq:boundary-orientation-compatibility}, \eqref{eq:immersed-eta-contraction} and then \eqref{eq:tau_partial-Sigma-and-contraction}, we have
\begin{align*}
\iota_{\eta(F(z),\omega_\Sigma(z))}\omega_\Sigma(z)
=(-1)^{m-1}\tau_\gamma(F(z))
=(-1)^{m-1}\tau_{\partial\Sigma(z)}(z)
=\iota_{\eta_{\partial\Sigma}(z)}\omega_\Sigma(z),
\end{align*}
which implies that
\begin{align}\label{eq:eta-F-eta-gamma}
\eta_{\partial\Sigma}(z)
=\eta(F(z),\omega_\Sigma(z)).
\end{align}

The oriented Gauss map is thus $\omega_\Sigma:\Sigma\to\widetilde{\mathcal G}_{m,n}$,
we also put for simplicity $P(z):=\mathbf q(\omega_\Sigma(z))$. For $z\in{\rm int} \Sigma$ and $v\in\mathbb R^n$, we define the curvature map by
\begin{align}\label{eq:source-curvature-field}
\widetilde{B}(z)[v]
:=D_{P(z)v}\omega_\Sigma(z).
\end{align}
Clearly, one has
\begin{align}\label{eq:shape-operator-tangential}
\widetilde{B}(z)\in {\rm Hom}(\mathbb R^n,T_{\omega_\Sigma(z)}\widetilde{\mathcal G}_{m,n}).
\end{align}
We also set the bilinear form
\begin{align}\label{defn:B-unoriented}
B(z)
:=D\overline{\mathbf q}_\omega\circ\widetilde B(z)
\end{align}
where $\overline{\mathbf q}$ can be taken as any $C^1$ extension of $\mathbf q$ to a neighborhood of $\widetilde{\mathcal G}_{m,n}$ thanks to \eqref{eq:shape-operator-tangential}. The trace of $B$ is then denoted as $\mathcal T(\widetilde{B}):={\rm tr}B=\sum_{i=1}^n B(e_i,e_i)$.

To compute the norm of $\widetilde{B}$, we choose a local oriented orthonormal tangent frame $X=(\tau_1,\ldots,\tau_m)$ and put for each $i\in\{1,\cdots,m\}$
\begin{align*}
Z_i
:=\left(A(\tau_i,\tau_1),\ldots,A(\tau_i,\tau_m)\right).
\end{align*}
Note that $X^TX$ is the orthogonal projection onto $T\Sigma$, and hence by Gauss formula we have $A(\tau_i,\tau_j)=(I_n-XX^T)D_{\tau_i}\tau_j$, so that $Z_i$ can be rewritten as
\begin{align*}
Z_i
=(I_n-XX^T)D_{\tau_i}X,
\end{align*}
which implies that $X^TZ_i=0$, thus $Z_i\in\mathcal H_X$.
By the chain rule, we find
\begin{align*}
\widetilde B[\tau_i]
=D_{\tau_i}\omega_\Sigma
=d\widetilde\iota_{[X]_{SO(m)}}(Z_i).
\end{align*}
Since $\widetilde\iota$ is isometric and $\widetilde{B}$ vanishes on normal vectors by definition, we thus obtain
\begin{align}\label{eq:widetilde-B=A-len}
|\widetilde{B}|^2
&=\sum_{i=1}^m g_{[X]_{SO(m)}}(Z_i,Z_i)\\
&=\sum_{i=1}^m{\rm tr}(Z_i^TZ_i)
=\sum_{i,j=1}^m|A(\tau_i,\tau_j)|^2
=|A|^2.
\end{align}

For any test vector field $\phi\in C_c^1(U\times\Lambda_m(\mathbb R^n);\mathbb R^n)$, we define $X(z):=\phi(F(z),\omega_\Sigma(z))$ and $X^T(z):=P(z)X(z)$. By the chain rule we get
\begin{align}\label{eq:smooth-cv-ptws}
{\rm div} X^T
=D_\omega\phi\cdot \widetilde{B}+\left<\mathcal T(\widetilde{B}),\phi\right>+\left< D_x\phi,P\right>,
\end{align}
where all terms involving $\phi$ are evaluated at $(F(z),\omega_\Sigma(z))$, and we adopt the convention $D_\omega\phi\cdot \widetilde{B}:={\rm tr}(D_\omega\phi\circ \widetilde{B})$ as in \cite{KM22}.
By an approximation argument as in \cite[Section 2]{KM22}, we get by using the tangential divergence theorem that
\begin{align}\label{eq:tangential-divergence-F}
\int_\Sigma{\rm div} X^T \,d\mathcal H^m
=-\int_{\partial\Sigma}\left<\phi(F(z),\omega_\Sigma(z)),\eta_{\partial\Sigma}(z)\right>d\mathcal H^{m-1}.
\end{align}

Using the map $\Phi(z):=(F(z),\omega_\Sigma(z))$, we can define the oriented varifolds associated to $\Sigma$ and $\partial\Sigma$ as $V:=\Phi_\#(\mathcal H^m\llcorner\Sigma)$, $\Gamma:=(\Phi|_{\partial\Sigma})_\#(\mathcal H^{m-1}\llcorner\partial\Sigma)$.
Combining \eqref{eq:tangential-divergence-F}, \eqref{eq:smooth-cv-ptws} and \eqref{eq:eta-F-eta-gamma}, we obtain
\begin{align}
\int_{U\times\widetilde{\mathcal G}_{m,n}}
\left(D_\omega\phi\cdot\widetilde{B}+\left<\mathcal T(\widetilde{B}),\phi\right>
+\left< D_x\phi,\mathbf q(\omega)\right>\right) dV
=-\int_{\widetilde G_{m,n,\gamma}}\left<\phi,\eta\right> d\Gamma.
\end{align}
This motivates the following definition.

\begin{definition}\label{defn:ocv}
\normalfont
Let $V$ be an oriented $m$-varifold in $\mathbb{R}^n$, let $\mathcal N$ be a tubular neighborhood of $\widetilde{\mathcal G}_{m,n}$ in $\Lambda_m(\mathbb R^n)$, and let $\overline{\mathbf q}:\mathcal N\to{\rm Sym}_n(\mathbb R)$ be a $C^1$-extension of $\mathbf q$ to $\mathcal N$. We say that $V$ is an {\em oriented curvature $m$-varifold with
given oriented boundary $\gamma$}, if there exist $\widetilde{B}\in L_{\mathrm{loc}}^1\left(V;{\rm Hom}(\mathbb R^n,\Lambda_m(\mathbb R^n))\right)$ and a Radon measure $\Gamma$ on $\widetilde G_{m,n,\gamma}$ such that for any test function $\phi\in C_c^1\left(\mathbb R^n\times\Lambda_m(\mathbb R^n);\mathbb R^n\right)$, the curvature identity holds:
\begin{align}\label{iden:ocv}
\int_{\mathbb R^n\times\widetilde{\mathcal G}_{m,n}}
\left(D_\omega\phi\cdot \widetilde{B}+\left<\mathcal T(\widetilde{B}),\phi\right>
+\left< D_x\phi,\mathbf q(\omega)\right>\right) dV
=-\int_{\widetilde G_{m,n,\gamma}}\left<\phi,\eta\right> d\Gamma.
\end{align}
\end{definition}

\begin{remark}
\normalfont
If $\Gamma=0$, one can show by a direct computation that the oriented curvature identity \eqref{iden:ocv} is consistent with the ones derived by Delladio-Scianna \cite{DS95} and Bellettini \cite{Bellettini23}.
\end{remark}

{
As an easy consequence of Definition \ref{defn:ocv}, we have the following {\em tangential property of curvature}.

\begin{proposition}\label{Prop:tangentiality-S}
Let $(V,\widetilde{B},\Gamma)$ satisfy Definition \ref{defn:ocv}. Then one has
\begin{align}\label{eq:S_i-tang}
\widetilde{B}(x,\omega)\in{\rm Hom}(\mathbb R^n, T_\omega\widetilde{\mathcal G}_{m,n})
\end{align}
for $V$-a.e. $(x,\omega)$. Consequently, Definition
\ref{defn:ocv} is independent of the chosen $C^1$-extension $\overline{\mathbf q}$ of $\mathbf q$, and
\begin{align}\label{eq:B-q_omega}
B(x,\omega)
:=D\overline{\mathbf q}_\omega\circ\widetilde B(x,\omega)
=D\mathbf q_\omega\circ\widetilde B(x,\omega)
\end{align}
\end{proposition}

\begin{proof}
Let $\psi\in
C_c^1\left(\mathbb R^n\times\Lambda_m(\mathbb R^n)\right)$ satisfy $\psi=0$ on $\mathbb R^n\times\widetilde{\mathcal G}_{m,n}$, so that for every $\omega\in\widetilde{\mathcal G}_{m,n}$, the function
$x\mapsto\psi(x,\omega)$ vanishes identically on $\mathbb R^n$, which implies: for every $(x,\omega)\in
\mathbb R^n\times\widetilde{\mathcal G}_{m,n}$,
\begin{align}\label{eq:D-psi}
D_x\psi(x,\omega)=0,\quad
D_\omega\psi(x,\omega)[\zeta]=0\quad\forall\zeta\in T_\omega\widetilde{\mathcal G}_{m,n}.
\end{align}
Fix any $e\in\mathbb R^n$. Testing \eqref{iden:ocv} with $\phi=\psi e$, and noting that 
\begin{align*}
D_\omega\phi\cdot\widetilde B
={\rm tr}(D_\omega\phi\circ \widetilde{B})
=\sum_{j=1}^n\left<D_\omega\phi[\widetilde B[e_j]]],e_j\right>
=\sum_{j=1}^nD_\omega\psi[\widetilde B[e_j]]\left<e,e_j\right>
=D_\omega\psi[\widetilde B[e]],
\end{align*}
we thus find
\begin{align}\label{eq:test-vani}
\int_{\mathbb R^n\times\widetilde{\mathcal G}_{m,n}}D_\omega\psi(x,\omega)[\widetilde{B}(x,\omega)[e]]
\,dV(x,\omega)
=0.
\end{align}
If we let $N_\omega\widetilde{\mathcal G}_{m,n}:=\left(T_\omega\widetilde{\mathcal G}_{m,n}
\right)^\perp\subset\Lambda_m(\mathbb R^n)$
and $\widetilde{B}^\perp(x,\omega):={\rm proj}_{N_\omega\widetilde{\mathcal G}_{m,n}}\circ\widetilde{B}(x,\omega)$,
then the integrand can be further written as
\begin{align*}
D_\omega\psi(x,\omega)[\widetilde{B}(x,\omega)[v]]
=D_\omega\psi(x,\omega)[\widetilde{B}^\perp(x,\omega)[v]].
\end{align*}

To proceed, let $Z:\mathbb R^n\times\widetilde{\mathcal G}_{m,n}\to\Lambda_m(\mathbb R^n)$ be a compactly supported $C^1$ map satisfying $Z(x,\omega)\in N_\omega\widetilde{\mathcal G}_{m,n}$ for every $(x,\omega)$.
Let $\mathfrak p:\mathcal N\to\widetilde{\mathcal G}_{m,n}$
be the nearest point projection in the tubular neighborhood $\mathcal N$, and define
\begin{align*}
\nu(\xi)
:=\xi-\mathfrak p(\xi),\quad\forall\xi\in\mathcal N,
\end{align*}
which satisfies $\nu\mid_{\widetilde{\mathcal{G}}_{m,n}}=0$.
We take a cut-off function $\chi\in C_c^1(\mathbb R^n\times\mathcal N)$ such that $\chi=1$ on a neighborhood of ${\rm spt}Z$, and put
\begin{align*}
\psi(x,\xi)
:=\chi(x,\xi)\left<\nu(\xi),Z\left(x,\mathfrak p(\xi)\right) \right>.
\end{align*}
After extending it by zero outside $\mathbb R^n\times\mathcal N$, the function $\psi$ belongs to $C_c^1\left(\mathbb R^n\times\Lambda_m(\mathbb R^n)\right)$ and satisfies $\psi=0$ on $\mathbb R^n\times\widetilde{\mathcal G}_{m,n}$.

At every $\omega\in\widetilde{\mathcal G}_{m,n}$, one clearly has
\begin{align*}
D\mathfrak p_\omega
={\rm proj}_{T_\omega\widetilde{\mathcal G}_{m,n}},\quad
D\nu_\omega
={\rm proj}_{N_\omega\widetilde{\mathcal G}_{m,n}}.
\end{align*}
Since $\nu(\omega)=0$, it follows that for every $(x,\omega)\in \mathbb R^n\times\widetilde{\mathcal G}_{m,n}$ and every $\zeta\in\Lambda_m(\mathbb R^n)$,
\begin{align*}
D_\xi\psi(x,\omega)[\zeta]
=\left<\zeta^\perp,Z(x,\omega)\right>
\end{align*}
where $\zeta^\perp:={\rm proj}_{N_\omega\widetilde{\mathcal G}_{m,n}}\zeta$.
Applying \eqref{eq:test-vani} to this choice of
$\psi$, we obtain
\begin{align*}
\int_{\mathbb R^n\times\widetilde{\mathcal G}_{m,n}}\left<\widetilde B^\perp(x,\omega)[v],Z(x,\omega)\right> dV(x,\omega)
=0.
\end{align*}
Since $v$ is arbitrary, we deduce that $\widetilde B^\perp=0$ for $V$-a.e. $(x,\omega)$, which proves \eqref{eq:S_i-tang}.

Finally, let $\overline{\mathbf q}_1$ and $\overline{\mathbf q}_2$ be two $C^1$-extensions of $\mathbf q$ to neighborhoods of $\widetilde{\mathcal G}_{m,n}$, since their restrictions to $\widetilde{\mathcal G}_{m,n}$ agree, we clearly have $D\overline{\mathbf q}_{1,\omega}[\zeta]=D\overline{\mathbf q}_{2,\omega}[\zeta]=D\mathbf q_\omega[\zeta]$ for every
$\zeta\in T_\omega\widetilde{\mathcal G}_{m,n}$.
Taking \eqref{eq:S_i-tang} into account, we see that \eqref{eq:B-q_omega} holds, and Definition \ref{defn:ocv} (in particular, the curvature identity \eqref{iden:ocv}) is independent of the extension $\overline{\mathbf q}$.
\end{proof}

The proposition also shows that the curvature identity \eqref{iden:ocv} depends only on the restriction of the ambient test field to the oriented Grassmannian $\mathbb R^n\times\widetilde{\mathcal G}_{m,n}$.

\begin{corollary}\label{Cor:exten-indep}
Under the assumptions of Definition \ref{defn:ocv}, both sides of \eqref{iden:ocv} depend only on the restriction of $\phi\in C_c^1\left(\mathbb R^n\times\Lambda_m(\mathbb R^n);\mathbb R^n\right)$ to $\mathbb R^n\times\widetilde{\mathcal G}_{m,n}$.
\end{corollary}

\begin{proof}
Let $\phi_1,\phi_2\in C_c^1\left(\mathbb R^n\times\Lambda_m(\mathbb R^n);
\mathbb R^n\right)$ satisfy $\phi_1=\phi_2$ on $\mathbb R^n\times\widetilde{\mathcal G}_{m,n}$, and set $\psi:=\phi_1-\phi_2$, then $\psi=0$ on $\mathbb R^n\times\widetilde{\mathcal G}_{m,n}$.
By \eqref{eq:D-psi} and \eqref{eq:S_i-tang}, we find
\begin{align*}
D_\omega\psi(x,\omega)\circ\widetilde B(x,\omega)
=0\quad V\text{-a.e.}
\end{align*}
Back to \eqref{iden:ocv}, we see that each term on both sides vanishes, which completes the proof.
\end{proof}
}

{
\subsection{Basic properties}
\label{Sec:3-basic-prop}

We collect some direct consequences of Definition
\ref{defn:ocv}. Throughout this subsection,
let $(V,\widetilde{B},\Gamma)$ satisfy Definition \ref{defn:ocv}.
We recall that
\begin{align}\label{defn:trace-B}
\mathcal T(\widetilde{B})(x,\omega)={\rm tr}B(x,\omega)
\end{align}
with $B$ given by \eqref{eq:B-q_omega}.
Since $\mathbf q$ is smooth and $\widetilde{\mathcal G}_{m,n}$ is compact, there exists $C=C(m,n)>0$ such that
\begin{align}\label{eq:T-S-bound}
|\mathcal T(\widetilde B)(x,\omega)|
\leq C|\widetilde B(x,\omega)|,\quad V\text{-a.e }(x,\omega).
\end{align}
We observe the following simple facts:
\begin{enumerate}
\item Let $\pi:\mathbb R^n\times\Lambda_m(\mathbb R^n)\to U$ denote the canonical projection onto the first variable, and put $\|\Gamma\|:=\pi_\#\Gamma$.
By Lemma \ref{lem:disintegration}, there exists a $\|\Gamma\|$-measurable family of Radon probability measures $\left\{\Gamma^x\right\}_{x\in\gamma}$, such that $\Gamma^x$ is supported on $\widetilde G_{m,n,\gamma}(x)$ for $\|\Gamma\|$-a.e. $x$, and
\begin{align}
\int_{\widetilde G_{m,n,\gamma}}f(x,\omega)\,d\Gamma(x,\omega)
=\int_\gamma\int_{\widetilde G_{m,n,\gamma}(x)}
f(x,\omega)\,d\Gamma^x(\omega)d\|\Gamma\|(x)
\end{align}
for every nonnegative Borel function $f$.

We define the {\em generalized inward-pointing conormal} associated with $\Gamma$ by
\begin{align}\label{defn:eta_V}
\eta_V(x)
:=\int_{\widetilde G_{m,n,\gamma}(x)}\eta(x,\omega)\,d\Gamma^x(\omega)
\end{align}
for $\|\Gamma\|$-a.e. $x\in\gamma$. Then $\eta_V$ is $\|\Gamma\|$-measurable and satisfies
\begin{align*}
|\eta_V(x)|
\leq 1\quad\text{for }\|\Gamma\|\text{-a.e. }x.
\end{align*}
Note that $\eta_V$ is in general not unit, since different oriented tangent planes over the same boundary point may cause cancellation in the average.

\item Put $\|V\|:=\pi_\#V$. For any $\varphi\in C_c^1(\mathbb R^n;\mathbb R^n)$, testing \eqref{iden:ocv} with $\phi(x,\omega)=\varphi(x)$, we obtain
\begin{align}\label{eq:1st-variation-formu}
\delta V^0(\varphi)
&=-\int_{\mathbb R^n\times\widetilde{\mathcal G}_{m,n}}\left<\mathcal T(\widetilde B)(x,\omega),\varphi(x) \right> dV(x,\omega)\\
&\quad-\int_{\widetilde G_{m,n,\gamma}}\left<\varphi(x),\eta(x,\omega) \right> d\Gamma(x,\omega).
\end{align}
Here $\delta V^0(\varphi)=\int_{\mathbb R^n\times\mathcal P_{m,n}}
\left<D\varphi(x),P\right>dV^0(x,P)$ is the first variation of the unoriented varifold $V^0$.
By Lemma \ref{lem:disintegration}, there exists a $\|V\|$-measurable family of Radon probability measures $\{V^x\}_{x\in U}$ on $\widetilde{\mathcal G}_{m,n}$.
We can define
\begin{align}
H_{\widetilde B}(x)
:=\int_{\widetilde{\mathcal G}_{m,n}}\mathcal T(\widetilde B)(x,\omega)\,dV^x(\omega)
\end{align}
for $\|V\|$-a.e. $x$. By \eqref{eq:T-S-bound}, $H_{\widetilde B}$ is $\|V\|$-measurable and $H_{\widetilde B}\in L_{\rm loc}^1(\|V\|;\mathbb R^n)$. Using \eqref{defn:eta_V}, identity \eqref{eq:1st-variation-formu} then becomes
\begin{align}
\delta V^0(\varphi)
=-\int_\mathbb R^n\left<H_{\widetilde B}(x),\varphi(x) \right> d\|V\|(x)
-\int_\gamma\left<\eta_V(x),\varphi(x)\right> d\|\Gamma\|(x).
\end{align}
\end{enumerate}

We have the following scaling property, see also \cite{KM22}.
\begin{remark}\label{rem:scaling-prop}
\normalfont
Let $(V,\widetilde B,\Gamma)$ satisfy Definition \ref{defn:ocv}. For every $\lambda>0$, the triple $(\lambda V,\widetilde B,\lambda\Gamma)$ also satisfies Definition \ref{defn:ocv}. Namely, $\lambda V$ has the same weak curvature field $\widetilde B$ and boundary measure $\lambda\Gamma$.
Moreover, for every $1\leq p<\infty$ we have
\begin{align}
\|\widetilde B\|_{L^p(\lambda V)}
=\lambda^{1/p}\|\widetilde B\|_{L^p(V)}.
\end{align}
\end{remark}
It is also standard to observe the following fact.
\begin{proposition}
Let $(V,\widetilde B,\Gamma)$ satisfy Definition \ref{defn:ocv}. Then the unoriented varifold $V^0=\mathbf q_\#V$ has locally bounded first variation. Moreover, the density $\Theta^m(\|V\|,x)$ exists for $\|V\|$-a.e. $x\in \mathbb R^n$. If, in addition, there exists $a>0$ such that $\Theta^m(\|V\|,x)\geq a$ for $\|V\|$-a.e. $x$, then $V^0$ is a rectifiable $m$-varifold.
\end{proposition}

\begin{proof}
Let $K\subset\mathbb R^n$ be compact, let $\varphi\in C_c^1(\mathbb R^n;\mathbb R^n)$ be a test function satisfying $|\varphi|\leq1$ with ${\rm spt}\varphi\subset K$. By \eqref{eq:1st-variation-formu} and
\eqref{eq:T-S-bound},
\begin{align}\label{esti:1st-vari}
\left|\delta V^0(\varphi)\right|
&\leq\int_{K\times\widetilde{\mathcal G}_{m,n}}|\mathcal T(\widetilde B)(x,\omega)|
\,dV(x,\omega)+\Gamma\left(\pi^{-1}(K\cap\gamma)\right)\\
&\leq C\int_{K\times\widetilde{\mathcal G}_{m,n}}|\widetilde B(x,\omega)|\,dV(x,\omega)+\|\Gamma\|(K\cap\gamma),
\end{align}
where $C=C(m,n)>0$. The right-hand side is bounded since $\widetilde B\in L_{\rm loc}^1(V)$ and $\Gamma$ is a Radon measure. This shows that $V^0$ has locally bounded first variation.

Recall that $\|V^0\|=\|V\|$, the existence of
$\Theta^m(\|V\|,x)$ for $\|V\|$-a.e. $x$ therefore follows from \cite[Lemma 40.5]{Simon83}. Under the additional positive lower
density assumption, we can use \cite[Theorem 42.4]{Simon83} to deduce the rectifiability of $V^0$.
This completes the proof.
\end{proof}
}

\begin{proposition}\label{Prop:bdry-spt}
Let $(V,\widetilde B,\Gamma)$ satisfy Definition \ref{defn:ocv}. Then one has
$${\rm spt}(\pi_\#\Gamma)\subset{\rm spt}\|V\|.$$
\end{proposition}
\begin{proof}
We write for simplicity the open set $U_0:=\mathbb R^n\setminus{\rm spt}\|V\|$.
Testing the curvature identity \eqref{iden:ocv} with any $\phi\in C_c^1(U_0\times\Lambda_m(\mathbb R^n))$, we obtain
\begin{align*}
\int\left<\phi,\eta\right>d\Gamma
=0,
\end{align*}
which shows that the vector-valued Radon measure $\eta\Gamma$ vanishes on $U_0\times\Lambda_m(\mathbb R^n)$.
Since $\eta$ is unit on $\widetilde G_{m,n,\gamma}$, it then follows that $\pi_\#\Gamma(U_0)=0$, and the proof is completed.
\end{proof}

{
\subsection{Relation with unoriented curvature varifolds}\label{Sec:3-relation-unoriented}

In \cite{Hutchinson86}, Hutchinson introduced the notion of weak curvature for unoriented integral varifolds. This notion was later extended by Mantegazza \cite{Mantegazza96} to rectifiable curvature varifolds with boundary. The corresponding weak curvature identity also makes sense for general varifolds, without assuming rectifiability in advance. We now show that the varifolds in Definition
\ref{defn:ocv} induce the corresponding
unoriented curvature varifolds after forgetting the orientation (namely, after push-forward by $\mathbf{q}$).

Let $(V,\widetilde B,\Gamma)$ satisfy Definition
\ref{defn:ocv}. We decompose as in Section \ref{Sec:oriented-structure}:
\begin{align*}
V=V^s+V^{as},\quad
\Gamma=\Gamma^s+\Gamma^{as},
\end{align*}
such that $V^{as}=\sigma_VV^s$ and $\Gamma^{as}=\sigma_\Gamma\Gamma^s$, where $\sigma_V$ and $\sigma_\Gamma$ are functions resulting from Proposition \ref{Prop:sigma-odd}.
We define a vector-valued Radon measure
$\mathcal M$ on $\mathbb R^n\times\widetilde{\mathcal G}_{m,n}$ by letting
\begin{align}\label{defn:mcM}
\int\left<\Psi(x,\omega),d\mathcal M(x,\omega) \right>
:=\int\left<\Psi(x,\omega),D\mathbf q_\omega\circ\widetilde B(x,\omega) \right> dV(x,\omega)
\end{align}
for every $\Psi\in C_c(\mathbb R^n\times\Lambda_m(\mathbb R^n);\mathbb R^{n\times n})$.
The symmetric and anti-symmetric parts of $\mathcal M$ are then given by
\begin{align}\label{eq:mcM-decomp}
\mathcal M^s
:=\frac12\left(\mathcal M+\mathcal I_\#\mathcal M\right),\quad
\mathcal M^{as}
:=\frac12\left(\mathcal M-\mathcal I_\#\mathcal M\right).
\end{align}
Since $\mathbf q\circ\mathcal I=\mathbf q$, one has $\mathbf q_\#\mathcal M=\mathbf q_\#\mathcal M^s$, $\mathbf q_\#\mathcal M^{as}=0$.

Now we define $A\in L_{\mathrm{loc}}^1\left(V^0;{\rm Hom}(\mathbb R^n,{\rm Sym}_n(\mathbb R))\right)$ by
\begin{align}\label{defn:AV^0}
AV^0:=\mathbf q_\#(BV),
\end{align}
where $B$ is as in \eqref{eq:B-q_omega}.
Thanks to Proposition \ref{Prop:tangentiality-S} we have
\begin{align*}
A(x,P)\in{\rm Hom}(\mathbb R^n,T_P\mathcal P_{m,n})\text{ for }V^0\text{-a.e. }(x,P).
\end{align*}
On the other hand, since $\eta$ satisfies $\eta(x,-\omega)=-\eta(x,\omega)$ by \eqref{defn:eta(x,om)}, and the function $\sigma_\Gamma$ is also odd in $\omega$ (recalling Proposition \ref{Prop:sigma-odd}), we see that the vector-valued function $\sigma_\Gamma(x,\omega)\eta(x,\omega)$
is even in $\omega$. It therefore descends to a
$\Gamma^0$-measurable function, say $\widetilde\eta$, given by
\begin{align}\label{defn:unoriented-conorm}
\widetilde\eta\left(x,\mathbf q(\omega)\right)
:=\sigma_\Gamma(x,\omega)\eta(x,\omega).
\end{align}
In particular, $|\widetilde\eta|\leq1$ for $\Gamma^0$-a.e.

\begin{proposition}\label{prop:compat-Mantegazza}
Let $(V,\widetilde{B},\Gamma)$ satisfy Definition \ref{defn:ocv}. Then for every $\varphi\in C_c^1(\mathbb R^n\times\mathbb R^{n\times n};\mathbb R^n)$, the associated unoriented varifolds $V^0=\mathbf q_\#V$ and $\Gamma^0=\mathbf q_\#\Gamma$ satisfy
\begin{align}\label{iden:unoriented-CV}
\int_{\mathbb R^n\times\mathcal P_{m,n}}\left(D_P\varphi\cdot A+\left<{\rm tr}A,\varphi\right>+\left<D_x\varphi,P\right>\right)dV^0
=-\int_{G_{m,n,\gamma}}\left<\varphi,\widetilde\eta\right>d\Gamma^0,
\end{align}
where $G_{m,n,\gamma}:=\mathbf q(\widetilde G_{m,n,\gamma})$ and $D_P\varphi\cdot A:={\rm tr}(D_P\varphi\circ A)$. This agrees with the curvature identity in \cite{Mantegazza96}, see also \cite{KM22}.
\end{proposition}
\begin{proof}
Let $\varphi\in C_c^1\left(\mathbb R^n
\times\mathbb R^{n\times n};\mathbb R^n\right)$, and define the even test function $\phi(x,\omega):=\varphi\left(x,\mathbf q(\omega)\right)$. By the chain rule, we have $D_x\phi=D_x\varphi$ and also $D_\omega\phi\cdot\widetilde{B}=D_P\varphi\cdot B$, where the terms involving $\varphi$ are evaluated at $(x,\mathbf q(\omega))$.
Substituting this choice of $\phi$ into \eqref{iden:ocv} and using \eqref{defn:AV^0}, the interior terms thus read
\begin{align*}
\int_{\mathbb R^n\times\mathcal P_{m,n}}
\left(D_P\varphi\cdot A+\left<{\rm tr}A,\varphi\right>+\left< D_x\varphi,P\right>\right) dV^0.
\end{align*}
It remains to consider the boundary term. Since $\Gamma=(1+\sigma_\Gamma)\Gamma^s$, we have
\begin{align}
&\int_{\widetilde G_{m,n,\gamma}}\left<\varphi\left(x,\mathbf q(\omega)\right),\eta(x,\omega)\right>d\Gamma(x,\omega)\\
=&\int_{\widetilde G_{m,n,\gamma}}\left<\varphi\left(x,\mathbf q(\omega)\right),\eta(x,\omega)\right>\left(1+\sigma_\Gamma(x,\omega)\right)d\Gamma^s(x,\omega).
\end{align}
Since $\Gamma^s$ is invariant under $\mathcal I$, while $\eta(x,\omega)$ is odd and
$\varphi(x,\mathbf q(\omega))$ is even, we thus find
\begin{align}
&\int_{\widetilde G_{m,n,\gamma}}\left<\varphi\left(x,\mathbf q(\omega)\right),\eta(x,\omega) \right>
d\Gamma(x,\omega)\\
=&\int_{\widetilde G_{m,n,\gamma}}\left<\varphi\left(x,\mathbf q(\omega)\right),\sigma_\Gamma(x,\omega)\eta(x,\omega) \right> d\Gamma^s(x,\omega)\\
=&\int_{G_{m,n,\gamma}}\left<\varphi(x,P),\widetilde\eta(x,P) \right> d\Gamma^0(x,P).
\end{align}
Combining the above computations, we obtain \eqref{iden:unoriented-CV} as desired.
\end{proof}

\begin{remark}\label{rem:odd-test-funcs}
\normalfont
In the proof of Proposition \ref{prop:compat-Mantegazza} we only use test functions that are even in $\omega$, which amounts to be forgetting about the orientation in the curvature identity \eqref{iden:ocv}.
We note, however, that \eqref{iden:ocv} in fact contains more information, since one can test it with test functions which are odd in $\omega$. We delay the detailed discussion to the next subsection.
\end{remark}

}

\subsection{Odd test functions}

We continue to use the notations in Section \ref{Sec:3-relation-unoriented}.
As discussed in Remark \ref{rem:odd-test-funcs}, we wish to test the curvature identity \eqref{iden:ocv} with test functions $\phi(x,\omega)$ that are odd in $\omega$.
To this end, we let $\mathcal O\subset\mathcal P_{m,n}$ be any open set on which there exists a smooth lift of $\mathbf q$, say $\ell:\mathcal O\to\widetilde{\mathcal G}_{m,n}$, satisfying
\begin{align*}
\mathbf q(\ell(P))
=P,\quad\forall P\in\mathcal O.
\end{align*}
Moreover, we define a function $s:\ell(\mathcal O)\to\{-1,1\}$ by letting
\begin{align}\label{defn:s}
s(\ell(P))
=1,\quad
s(-\ell(P))
=-1.
\end{align}

For $A'$ given by $A'V^0:=\mathbf q_\#(s\mathcal M)$ where $\mathcal M$ is defined as \eqref{defn:mcM}, we note that by \eqref{eq:mcM-decomp} and that $s$ is odd on $\ell(\mathcal O)$, one has $A'\in L^1_{loc}(V^0;{\rm Hom}(\mathbb R^n,T_P\mathcal P_{m,n}))$ and that
\begin{align*}
A'V^0
:=\mathbf q_\#(s\mathcal M)
=\mathbf q_\#(s\mathcal M^{as}).
\end{align*}

\begin{proposition}
For any $\varphi\in C^1_c(\mathbb R^n\times\mathcal O;\mathbb R^n)$, one has
\begin{align}\label{iden:odd-CV}
&\int
D_P\varphi\cdot A'+\left<{\rm tr}A',\varphi\right>+\sigma_V(x,\ell(P))\left<D_x\varphi,P\right>dV^0(x,P)\\
=&-\int_{G_{m,n,\gamma}}\left<\varphi,\eta(x,\ell(P))\right>d\Gamma^0(x,P).
\end{align}
\end{proposition}
\begin{proof}
Let $\phi(x,\omega)\in C^1_c(\mathbb R^n\times\Lambda_m(\mathbb R^n);\mathbb R^n)$ be an odd test function, which agrees with $s(\omega)\varphi(x,\mathbf q(\omega))$ on $\mathbb R^n\times\mathbf q^{-1}(\mathcal O)$, and $\equiv0$ otherwise. Thanks to \eqref{defn:s}, $s$ is constant on the connected component of $\mathbf q^{-1}(\mathcal O)$, and hence by \eqref{eq:B-q_omega},
\begin{align*}
D_\omega\phi\cdot\widetilde B
=sD_P\varphi\cdot B.
\end{align*}
Testing the curvature identity \eqref{iden:ocv} with such a choice of $\phi$, we find
\begin{align}\label{iden:ocv-0}
\int sD_P\varphi\cdot B+s\left<{\rm tr}B,\varphi\right>+s\left<D_x\varphi,\mathbf q(\omega)\right>dV
=-\int s\left<\varphi,\eta\right>d\Gamma.
\end{align}
Since $s$ is odd on $\omega$, we have $\mathbf q_\#(sV)=\mathbf q_\#(sV^{as})=\mathbf q_\#(s\sigma_VV^s)$. Note that $s\sigma_V$ is even on $\omega$, hence it descends to $\sigma_V(x,\ell(P))$ for $\mathbf q_\#V^s$-a.e., thanks to \eqref{defn:s}.
Similarly, since $\eta$ is odd on $\omega$, we see that $\mathbf q_\#(s\eta\Gamma)=\mathbf q_\#(s\eta\Gamma^s)$, where $s\eta$ is even on $\omega$ and descends to $\eta(x,\ell(P))$ for $\mathbf q_\#\Gamma^s$-a.e., thanks again to \eqref{defn:s}.

Putting these facts back into \eqref{iden:ocv-0} and by definition of push-forward of measures, we obtain \eqref{iden:odd-CV}.
The proof is thus completed.
\end{proof}

It is well-known that oriented curvature varifolds are stronger than unoriented curvature varifolds, since there are many examples which satisfy the definitions in \cite{Hutchinson86,Mantegazza96} but not satisfy Definition \ref{defn:ocv}, see e.g. \cite[Remark 4.6]{DS95}.

These examples fail to satisfy Definition \ref{defn:ocv} but still satisfy \cite{Hutchinson86,Mantegazza96} since a sudden change of the orientation of varifolds will violate the oriented curvature identity \eqref{iden:odd-CV}, but still satisfy the unoriented one \eqref{iden:unoriented-CV}.
To see this, consider any $(V,\widetilde{B},\Gamma)$ satisfying Definition \ref{defn:ocv},
if $\widetilde B=0$ for $V$-a.e., then for any $\varphi$ with $\varphi|_\gamma=0$, \eqref{iden:odd-CV} reduces to
\begin{align}
\int\sigma_V(x,\ell(P))\left<D_x\varphi,P\right>dV^0(x,P)
=0.
\end{align}
In particular, suppose that $V^0=\theta\mathcal H^m\llcorner M_0\otimes\delta_{P_0}$ locally in $U_0\subset\mathbb R^n$, where $U_0\cap\gamma=\emptyset$, $M$ is a connected open subset of an affine $m$-plane $P_0$, and $\theta$ is locally a constant on $M_0$, then for any $\varphi$ supported in $U_0$, we have
\begin{align}
\int_{M_0}\left(1+\sigma_V(x,\ell(P_0))\right)\left<D_x\varphi,P_0\right>d\mathcal H^m(x)
=0,
\end{align}
where $\int_{M_0}\left<D_x\varphi,P_0\right>\,d\mathcal H^m(x)=0$ by divergence theorem. Applying the Constancy theorem (Lemma \ref{Lem:Constancy}) to the non-negative measure $(1+\sigma_V(x,\ell(P_0)))\mathcal H^m\llcorner M_0$, we conclude that $\sigma_V(x,\ell(P_0))$ is constant $\mathcal H^m$-a.e. on $M_0$, which implies that there cannot be a sudden change of orientation.

\section{The case of zero curvature}\label{sec:4}
{
Throughout the section, let $K$ be a compact set in $\mathbb R^n$. We begin with a disintegration procedure using Lemma \ref{lem:disintegration} for oriented varifolds, which generalizes the one used in \cite[Section 3]{KM22}. For $\omega\in\widetilde{\mathcal G}_{m,n}$, we put $\mathbb P_\omega:={\rm span}(\omega)$, and set
\begin{align*}
p_\omega
:=\mathbf q(\omega),
\quad
p_\omega^\perp
:=I-p_\omega.
\end{align*}
Namely, $p_\omega$ and $p_\omega^\perp$ are the orthogonal projections
onto $\mathbb P_\omega$ and $\mathbb P_\omega^\perp$, respectively.

Let $V$ be an oriented $m$-varifold on $K\subset\mathbb R^n$, by standard disintegration (Lemma \ref{lem:disintegration}) we have the following:
\begin{enumerate}
\item Let $p_2:K\times\Lambda_m(\mathbb R^n)\to\Lambda_m(\mathbb R^n)$ be the projection onto the second factor, set
\begin{align*}
\nu:=(p_2)_\#V.
\end{align*}
There exists a $\nu$-measurable family of probability measures $\{\mu_\omega\}_{\omega\in\widetilde{\mathcal G}_{m,n}}$ on $K$, such that for every $\phi\in C_c^0\left(\mathbb R^n\times\Lambda_m(\mathbb R^n)\right)$,
\begin{align}\label{formu-V-measures-0}
\int_{K\times\widetilde{\mathcal G}_{m,n}}\phi(x,\omega)\,dV(x,\omega)
=\int_{\widetilde{\mathcal G}_{m,n}}\int_K\phi(x,\omega)\,d\mu_\omega(x)d\nu(\omega)
\end{align}

\item For $\nu$-a.e. $\omega$, we define
\begin{align*}
\lambda_\omega
:=(p_\omega^\perp)_\#\mu_\omega,
\end{align*}
which is a probability measure on $\mathbb P_\omega^\perp$. By disintegration there exists a $\lambda_\omega$-measurable family of probability measures $\{\mu_{\omega,z}\}_{z\in\mathbb P_\omega^\perp}$ such that $\mu_{\omega,z}$ is supported on $(z+\mathbb P_\omega)\cap K$ for $\lambda_\omega$-a.e. $z$, and for every $\varphi\in C_c^0(\mathbb R^n)$,
\begin{align}\label{eq:2nd-disintegra}
\int_K\varphi(x)\,d\mu_\omega(x)
=\int_{\mathbb P_\omega^\perp}\int_{(z+\mathbb P_\omega)\cap K}\varphi(x)\,d\mu_{\omega,z}(x)d\lambda_\omega(z).
\end{align}

\item For every $\varphi\in C_c^0(\mathbb R^n)$ and $\psi\in C^0(\widetilde{\mathcal G}_{m,n})$, one has
\begin{align}\label{formu-V-measures}
&\int_{K\times\widetilde{\mathcal G}_{m,n}}
\varphi(x)\psi(\omega)\,dV(x,\omega)\\
=&\int_{\widetilde{\mathcal G}_{m,n}}\psi(\omega)\int_{\mathbb P_\omega^\perp}\int_{(z+\mathbb P_\omega)\cap K}\varphi(x)\,d\mu_{\omega,z}(x)d\lambda_\omega(z)d\nu(\omega).
\end{align}
\end{enumerate}

\subsection{Oriented curvature varifolds with zero curvature and given oriented boundary}

Let $\gamma\subset K$ be a compact, connected, oriented, embedded $C^1$ submanifold of dimension $m-1$ without boundary. We continue to use the notations above, and for $\omega\in\widetilde{\mathcal G}_{m,n}$, $z\in\mathbb P_\omega^\perp$, we write for simplicity
\begin{align}\label{defn:L_omega,z}
L_{\omega,z}
:=z+\mathbb P_\omega.
\end{align}

\begin{lemma}\label{lem:KM22-Lem3.2}
Let $(V,\widetilde B,\Gamma)$ satisfy Definition
\ref{defn:ocv}, with ${\rm spt}\|V\|\subset K$. Assume that $\widetilde B=0$ for $V$-a.e.
Then for $\nu$-a.e. $\omega\in\widetilde{\mathcal G}_{m,n}$ and
$\lambda_\omega$-a.e. $z\in\mathbb P_\omega^\perp$, one has
\begin{align}\label{eq:0-curvature-slice}
\int_{L_{\omega,z}}\left<D\varphi(x),p_\omega\right>d\mu_{\omega,z}(x)
=0
\end{align}
for every $\varphi\in C_c^1(\mathbb R^n;\mathbb R^n)$ satisfying either $\varphi=0$ on $\gamma$, or $\varphi(x)\in T_x\gamma$ for every $x\in\gamma$.
\end{lemma}

\begin{proof}
Let $\varphi\in C_c^1(\mathbb R^n;\mathbb R^n)$ satisfy one of the two conditions in the statement, and let $\psi\in C^1_c(\Lambda_m(\mathbb R^n))$. Testing \eqref{iden:ocv} with $\phi(x,\omega):=\varphi(x)\psi(\omega)$ and using $\widetilde B=0$, we get
\begin{align}
\int_{K\times\widetilde{\mathcal G}_{m,n}}\psi(\omega)\left<D\varphi(x),p_\omega\right> dV(x,\omega)
=-\int_{\widetilde G_{m,n,\gamma}}\psi(\omega)
\left<\varphi(x),\eta(x,\omega)\right> d\Gamma(x,\omega).
\end{align}
If $\varphi=0$ on $\gamma$, then the RHS vanishes. If $\varphi(x)\in T_x\gamma$ on $\gamma$, it also vanishes because $\eta(x,\omega)\perp T_x\gamma$ for every $(x,\omega)\in\widetilde G_{m,n,\gamma}$. Thus in both cases, we find
\begin{align*}
\int_{K\times\widetilde{\mathcal G}_{m,n}}\psi(\omega)\left<D\varphi(x),p_\omega\right>dV(x,\omega)
=0.
\end{align*}
Using \eqref{formu-V-measures-0}, we thus obtain
\begin{align*}
\int_{\widetilde{\mathcal G}_{m,n}}\psi(\omega)\int_K\left<D\varphi(x),p_\omega\right>d\mu_\omega(x)d\nu(\omega)
=0.
\end{align*}
Since $\psi$ is arbitrary, it follows that for $\nu$-a.e. $\omega$,
\begin{align}\label{eq:fixed-ome}
\int_K\left<D\varphi(x),p_\omega\right> d\mu_\omega(x)
=0.
\end{align}
As in \cite[Lemma 3.2]{KM22}, by using a countable dense family of test vector fields, we see that the null set of $\omega$ may be chosen independently of $\varphi$, and we fix such an $\omega$.

Let $h\in C_c^1(\mathbb P_\omega^\perp)$, we note that
\begin{align*}
D_\tau[h(p^\perp_\omega x)]
=Dh(p^\perp_\omega x)[p^\perp_\omega \tau]
=0,\quad\forall\tau\in\mathbb P_\omega,
\end{align*}
and hence
\begin{align*}
\left<D\bigl[h(p_\omega^\perp x)\varphi(x)\bigr],p_\omega\right>
=h(p_\omega^\perp x)\left<D\varphi(x),p_\omega\right>.
\end{align*}
Replacing $\varphi(x)$ in \eqref{eq:fixed-ome} by $h(p_\omega^\perp x)\varphi(x)$ then using \eqref{eq:2nd-disintegra}, we find
\begin{align*}
\int_{\mathbb P_\omega^\perp}h(z)\underbrace{\int_{L_{\omega,z}}\left<D\varphi(x),p_\omega\right>d\mu_{\omega,z}(x)}_{:=\mathcal I_1}d\lambda_\omega(z)
=0.
\end{align*}
Since $h$ is arbitrary, the inner integral $\mathcal I_1$ vanishes for $\lambda_\omega$-a.e. $z$. Arguing as above, we see that the null set of $z$ is independent of $\varphi$, which completes the proof.
\end{proof}

\begin{lemma}\label{lem:0-CV-Interior}
Under the assumptions of Lemma \ref{lem:KM22-Lem3.2}, for $\nu$-a.e. $\omega\in\widetilde{\mathcal G}_{m,n}$ and $\lambda_\omega$-a.e. $z\in\mathbb P_\omega^\perp$, the following holds. For every connected component $\mathscr U$ of $L_{\omega,z}\setminus\gamma$ with $L_{\omega,z}$ given by \eqref{defn:L_omega,z}, there exists a constant $\theta_{\mathscr U}\geq0$ such that
\begin{align}\label{eq:constancy}
\mu_{\omega,z}\llcorner\mathscr U
=\theta_{\mathscr U}\mathcal H^m\llcorner\mathscr U.
\end{align}
\end{lemma}

\begin{proof}
Let $\mathscr U$ be a connected component of $L_{\omega,z}\setminus\gamma$, and let $\hat e_1,\ldots,\hat e_m$ be an orthonormal basis of $\mathbb P_\omega$. For $f\in C_c^\infty(\mathscr U)$ and $a\in\{1,\ldots,m\}$, we extend the vector field $f\hat e_a$ from $L_{\omega,z}$ to a compactly supported $C^1$ vector field $\varphi$ on $\mathbb R^n$ whose support is disjoint from $\gamma$. By construction, we have
\begin{align*}
\left<D\varphi(x),p_\omega\right>
=D_{\hat e_a}f(x)\quad\text{ on }L_{\omega,z}.
\end{align*}
Back to \eqref{eq:0-curvature-slice}, we get
\begin{align*}
\int_{\mathscr U}D_{\hat e_a}f(x)\,d\mu_{\omega,z}(x)=0
\end{align*}
for every $a=1,\ldots,m$. The Constancy Theorem (Lemma \ref{Lem:Constancy}) then yields \eqref{eq:constancy} and completes the proof.
\end{proof}

\subsubsection{A dichotomy result}
The following topological dichotomy will be useful fur our purpose.

\begin{proposition}\label{Prop:-dichotomy}
Let $\gamma\subset\mathbb R^n$ be a compact, connected, oriented, embedded $C^1$ submanifold of dimension $m-1$ without boundary, and let $L\subset\mathbb R^n$ be an affine $m$-plane. Then exactly one of the following alternatives holds:
\begin{enumerate}
\item $\gamma\subset L$;
\item $L\setminus\gamma$ is connected and unbounded.
\end{enumerate}
In particular, if $\gamma\subset L$ then $L\setminus\gamma$ has exactly one bounded connected component and the other one unbounded.
\end{proposition}

\begin{proof}
If $\gamma\subset L$, then the assertion follows from the Jordan Brouwer theorem, so it suffices to prove (2). We assume that $\gamma\not\subset L$, then argue by contradiction and suppose that $L\setminus\gamma$ is disconnected.

First note that, since $\gamma$ is compact, there exists $R>0$ such that $\gamma\subset B_R(0)$. Now by $m\geq2$, the set $L\setminus\overline{B_R(0)}$ is connected, which implies that $L\setminus\gamma$ has a unique unbounded connected component. At the meantime, it is assumed to be disconnected, therefore it must have a nonempty {\em bounded} connected component, say $\Omega$.

\noindent{\em Claim. The relative boundary of $\Omega$ in $L$, say $\partial_L\Omega$, satisfies
\begin{align*}
\partial_L\Omega\subset\gamma\cap L.
\end{align*}
}
Indeed, if $x\in\partial_L\Omega$ and $x\notin\gamma$, then a sufficiently small neighborhood of $x$ in $L$ is contained in
$L\setminus\gamma$, namely, $x$ must be an interior point of a connected component of $L\setminus\gamma$, contradicting $x\in\partial_L\Omega$.
Since $\gamma$ is a compact $C^1$ submanifold of dimension $m-1$, we have $\mathcal H^{m-1}(\gamma)<\infty$. Thus $\mathcal H^{m-1}(\partial_L\Omega)<\infty$, and hence $\Omega$ is a set of finite perimeter in $L$.

To proceed, we identify $L$ with $\mathbb R^m$ and fix an orientation on $L$, then let $\llbracket\Omega\rrbracket$ be the $m$-dimensional current induced by $\Omega$, with the orientation induced from $L$.
Since $\Omega$ is a set of finite perimeter in $L$, we see that
$T:=\partial\llbracket\Omega\rrbracket$ is an integral $(m-1)$-current satisfying $\partial T=0$ and
\begin{align*}
{\rm spt}T
\subset\partial_L\Omega
\subset\gamma\cap L.
\end{align*}
Since $\gamma$ is connected, oriented, and has no boundary, the Constancy Theorem (cf. \cite[Theorem 26.27]{Simon83}) applied on $\gamma$ yields: there exists some $k\in\mathbb Z$, such that
\begin{align*}
T=k\llbracket\gamma\rrbracket.
\end{align*}
On the other hand, since ${\rm spt}T\subset L$,
if $k\neq0$, then we have ${\rm spt}T={\rm spt}\llbracket\gamma\rrbracket=\gamma$,
and hence $\gamma\subset L$, which contradicts the assumption. Therefore we must have $k=0$ and hence $T=\partial\llbracket\Omega\rrbracket=0$. That is to say, $D_L\chi_\Omega=0$. Since $L$ is connected, this implies that $\chi_\Omega$ is constant almost everywhere on $L$, which is impossible, since by construction, $\Omega$ is a nonempty bounded open subset of $L$.
Therefore, we conclude that $L\setminus\gamma$ is connected.

Finally, recall that $\gamma\cap L\subset B_R(0)\cap L$, and hence $L\setminus B_R(0)\subset L\setminus\gamma$, so that $L\setminus\gamma$ is unbounded. This shows the second alternative and completes the proof.
\end{proof}

\begin{definition}[$(m,\gamma)$-slice]\label{Defn:m-gamma-slice}
\normalfont
Let $L\subset\mathbb R^n$ be an affine $m$-plane. A connected component of $L\setminus\gamma$ is called an {\em $m$-slice} in $L$.

An $m$-slice $\Delta\subset L$ is called an {\em $(m,\gamma)$-slice} if $\Delta$ is bounded. We denote by $\mathscr S_{m,\gamma}(L)$ the set of $(m,\gamma)$-slices in $L$.
\end{definition}

By Proposition \ref{Prop:-dichotomy}, we have $\mathscr S_{m,\gamma}(L)=\emptyset$, unless $\gamma\subset L$. If $\gamma\subset L$, then $\mathscr S_{m,\gamma}(L)$ consists of the unique bounded connected component of $L\setminus\gamma$, whose boundary is $\gamma$.

\subsubsection{Structural theorem}
\begin{theorem}\label{Thm:structure-V-Gamma-S=0}
Let $\gamma\subset K\subset\mathbb R^n$ be a compact, connected, oriented, embedded $C^1$ submanifold of dimension $m-1$ without boundary. Let $(V,\widetilde B,\Gamma)$ satisfy Definition
\ref{defn:ocv}, with ${\rm spt}\|V\|\subset K$. Suppose that $\widetilde B=0$ for $V$-a.e. Then, for $\nu$-a.e. $\omega$ and $\lambda_\omega$-a.e. $z$, there exist constants $\theta_{\omega,z,\Delta}\geq0$ for any $\Delta\in\mathscr S_{m,\gamma}(L_{\omega,z})$, such that
\begin{align}\label{eq:measure-struc}
\mu_{\omega,z}
=\sum_{\Delta\in\mathscr S_{m,\gamma}(L_{\omega,z})}\theta_{\omega,z,\Delta}\mathcal H^m\llcorner\Delta,
\end{align}
and for any test function $\phi\in C_c^0\left(\mathbb R^n\times\Lambda_m(\mathbb R^n)\right)$, one has
\begin{align}\label{eq:V-struc}
V(\phi)
&=\int_{\widetilde{\mathcal G}_{m,n}}\int_{\mathbb P_\omega^\perp}\sum_{\Delta\in\mathscr S_{m,\gamma}(L_{\omega,z})}\theta_{\omega,z,\Delta}\int_\Delta\phi(x,\omega)\,d\mathcal H^m(x)d\lambda_\omega(z)d\nu(\omega).
\end{align}
\end{theorem}

\begin{proof}
We fix $\omega$ and $z$ outside the null sets in Lemma \ref{lem:0-CV-Interior}, it follows that that on any connected component
$\mathscr U$ of $L_{\omega,z}\setminus\gamma$, there exists $\theta_{\mathscr U}\geq0$ such that
\begin{align*}
\mu_{\omega,z}\llcorner\mathscr U
=\theta_{\mathscr U}
\mathcal H^m\llcorner\mathscr U.
\end{align*}

\noindent{\em Case 1. $\gamma\not\subset L_{\omega,z}$.}

By Proposition \ref{Prop:-dichotomy},
$L_{\omega,z}\setminus\gamma$ is connected and unbounded. Since $\gamma$ is compact, this component has infinite $\mathcal H^m$-measure. Since $\mu_{\omega,z}$ is a probability measure, its multiplicity on $L_{\omega,z}\setminus\gamma$ must vanish. Hence we conclude
\begin{align}
\mu_{\omega,z}\llcorner\left(L_{\omega,z}\setminus\gamma\right)
=0.
\end{align}

\noindent{\em Case 2. $\gamma\subset L_{\omega,z}$.}

In this case, thanks to Proposition \ref{Prop:-dichotomy} again, $L_{\omega,z}\setminus\gamma$ consists of a unique bounded component $\Delta\in\mathscr S_{m,\gamma}(L_{\omega,z})$, and a unique unbounded component. The unbounded component has infinite $\mathcal H^m$-measure, so the multiplicity of $\mu_{\omega,z}$ on it is zero. On $\Delta$, we write
\begin{align*}
\mu_{\omega,z}\llcorner\Delta
=\theta_{\omega,z,\Delta}\mathcal H^m\llcorner\Delta,
\end{align*}
and note that $\partial\Delta=\gamma$.

To proceed, we define $\widetilde V:=V\llcorner\left((\mathbb R^n\setminus\gamma)\times\widetilde{\mathcal G}_{m,n}\right)$. By \eqref{formu-V-measures-0}, \eqref{eq:2nd-disintegra} and the conclusions above, for every compactly supported continuous $\phi$, one has
\begin{align}\label{eq:tilde-V-slice}
\widetilde V(\phi)
&=\int_{\widetilde{\mathcal G}_{m,n}}\int_{\mathbb P_\omega^\perp}\sum_{\Delta\in\mathscr S_{m,\gamma}(L_{\omega,z})}\theta_{\omega,z,\Delta}\int_\Delta\phi(x,\omega)\,d\mathcal H^m(x)d\lambda_\omega(z)d\nu(\omega)
\end{align}

It remains to prove that $V$ has zero mass on $\gamma\times\widetilde{\mathcal G}_{m,n}$. To this end we set
\begin{align*}
\Lambda
:=V-\widetilde V
=V\llcorner\left(\gamma\times\widetilde{\mathcal G}_{m,n}\right).
\end{align*}
Clearly $\Lambda$ is a nonnegative Radon measure. Since $\gamma$ is $C^1$ and compact, we know that the mapping $x\to{\rm proj}_{N_x\gamma}$ is continuous on $\gamma$ with
\begin{align*}
|{\rm proj}_{N_x\gamma}(y-x)|
=o(|y-x|)
\end{align*}
uniformly for $x,y\in\gamma$. Thus we can use the $C^1$-Whitney extension theorem (cf. \cite[Section 15.2]{Mag12}) to obtain a function $\varphi\in C^1_c(\mathbb R^n;\mathbb R^n)$ with $\varphi=0$ on $\gamma$ and
\begin{align}\label{eq:C^1-Whitney-D-varphi}
D\varphi(x)
={\rm proj}_{N_x\gamma},\quad\forall x\in\gamma.
\end{align}
Now let $\chi\in C^1_c(\Lambda_m(\mathbb R^n))$ be a cutoff function which $\equiv1$ in a neighborhood of $\widetilde{\mathcal G}_{m,n}$. Since $\widetilde B=0$ and $\varphi=0$ on $\gamma$, testing the curvature identity \eqref{iden:ocv} with $\phi(x,\omega):=\chi(\omega)\varphi(x)$, we get
\begin{align}\label{eq:normal-test-V}
\int_{K\times\widetilde{\mathcal G}_{m,n}}
\left<D\varphi(x),p_\omega\right>dV(x,\omega)
=0.
\end{align}

On the other hand, using
\eqref{eq:tilde-V-slice} we obtain
\begin{align}\label{eq:normal-test-tilde-V}
&\int_{K\times\widetilde{\mathcal G}_{m,n}}\left<D\varphi(x),p_\omega\right>d\widetilde V(x,\omega)\\
=&\int_{\widetilde{\mathcal G}_{m,n}}\int_{\mathbb P_\omega^\perp}\sum_{\Delta\in\mathscr S_{m,\gamma}(L_{\omega,z})}\theta_{\omega,z,\Delta}\int_\Delta{\rm div}_{L_{\omega,z}}\varphi \,d\mathcal H^md\lambda_\omega(z)d\nu(\omega)\\
=&0,
\end{align}
where we have used the divergence theorem together with the facts that  $\partial\Delta=\gamma$ and $\varphi=0$ on $\gamma$ for the last equality.

Subtracting \eqref{eq:normal-test-tilde-V} from
\eqref{eq:normal-test-V}, we find
\begin{align*}
\int_{\gamma\times\widetilde{\mathcal G}_{m,n}}\left<D\varphi(x),p_\omega\right> d\Lambda(x,\omega)
=0.
\end{align*}
Note also that for $x\in\gamma$, by \eqref{eq:C^1-Whitney-D-varphi} we have
\begin{align}
\left<D\varphi(x),p_\omega\right>
={\rm tr}\left(p_\omega({\rm Id}-{\rm proj}_{T_x\gamma})\right)
=m-{\rm tr}(p_\omega {\rm proj}_{T_x\gamma})
\geq m-(m-1)=1.
\end{align}
Since $\Lambda$ is nonnegative, it follows that $\Lambda=0$, and hence $V=\widetilde V$, which proves \eqref{eq:measure-struc}.
\eqref{eq:V-struc} then follows from \eqref{formu-V-measures}. This finishes the proof.
\end{proof}

For every $\Delta\in\mathscr S_{m,\gamma}(L_{\omega,z})$, we let $\eta_\Delta^{\mathrm{in}}(x)$ denote the inward-pointing unit conormal to $\gamma=\partial\Delta$ within the affine plane $L_{\omega,z}$.

\begin{corollary}\label{Cor:structure-bdry}
Under the assumptions and notation of Theorem \ref{Thm:structure-V-Gamma-S=0}, the boundary measure $\Gamma$ is characterized as follows: for every $\phi\in C_c^0\left(\mathbb R^n\times\Lambda_m(\mathbb R^n)\right)$,
\begin{align}\label{formu:Gamma-struc}
\Gamma(\phi)
&=\int_{\widetilde{\mathcal G}_{m,n}}\int_{\mathbb P_\omega^\perp}\sum_{\Delta\in\mathscr S_{m,\gamma}(L_{\omega,z})}\theta_{\omega,z,\Delta}\int_\gamma\left<\eta_\Delta^{\mathrm{in}}(x),\eta(x,\omega)\right>\phi(x,\omega)\,d\mathcal H^{m-1}(x)d\lambda_\omega(z)d\nu(\omega).
\end{align}
\end{corollary}

\begin{proof}
Let $\phi\in C_c^0\left(\mathbb R^n\times\Lambda_m(\mathbb R^n)\right)$. Since $\widetilde G_{m,n,\gamma}$ is compact and $\eta$ is continuous on $\widetilde G_{m,n,\gamma}$, by a standard approximation argument there exists $\widetilde\phi_j\in C^1_c(\mathbb R^n\times\Lambda_m(\mathbb R^n);\mathbb R^n)$ such that
\begin{align}\label{defn:widetilde-phi_j}
\widetilde\phi_j\to\phi\eta\text{ uniformly on }\widetilde G_{m,n,\gamma}\text{ as }j\to\infty.
\end{align}
Using the curvature identity \eqref{iden:ocv}, $\widetilde B=0$, and also \eqref{eq:V-struc}, we find
\begin{align}
&\int_{\widetilde G_{m,n,\gamma}}\left<\widetilde\phi_j(x,\omega),
\eta(x,\omega)\right>d\Gamma(x,\omega)\\
=&-\int_{K\times\widetilde{\mathcal G}_{m,n}}\left<D_x\widetilde\phi_j(x,\omega),p_\omega\right> dV(x,\omega)\\
=&-\int_{\widetilde{\mathcal G}_{m,n}}\int_{\mathbb P_\omega^\perp}\sum_{\Delta\in\mathscr S_{m,\gamma}(L_{\omega,z})}\theta_{\omega,z,\Delta}\int_\Delta{\rm div}_{L_{\omega,z}}\widetilde\phi_j(x,\omega)\,d\mathcal H^m(x)d\lambda_\omega(z)d\nu(\omega)\\
=&\int_{\widetilde{\mathcal G}_{m,n}}\int_{\mathbb P_\omega^\perp}\sum_{\Delta\in\mathscr S_{m,\gamma}(L_{\omega,z})}\theta_{\omega,z,\Delta}\int_\gamma\left<\widetilde\phi_j(x,\omega),\eta_\Delta^{\mathrm{in}}(x) \right> d\mathcal H^{m-1}(x)d\lambda_\omega(z)d\nu(\omega).
\end{align}
Letting $j\to\infty$, by virtue of \eqref{defn:widetilde-phi_j} and that $\eta$ is unit, we conclude the proof.
\end{proof}
}

\section{Compactness and existence of minimizers}\label{sec:5}

{
Let $(V,\widetilde B,\Gamma)$ satisfy Definition \ref{defn:ocv}. For brevity, we write $
V\in\widetilde{\bf CV}_{m,\gamma}(\mathbb R^n)$, and $\mathbb M(V):=\|V\|(\mathbb R^n)$, $\mathbb M(\Gamma)=\|\Gamma\|(\mathbb R^n)$ for the total masses. For $1\leq p<\infty$, we define
\begin{align*}
\mathcal E_p(V,\widetilde B)
:=\int_{\mathbb R^n\times\widetilde{\mathcal G}_{m,n}}|\widetilde B(x,\omega)|^p\,dV(x,\omega).
\end{align*}

\subsection{Compactness}

\begin{lemma}\label{Lem:bdry-mass-estimate}
For every $p\in[1,\infty)$, there exists a positive constant $C=C(m,n,p,\gamma)$ such that every $
V\in\widetilde{\bf CV}_{m,\gamma}(\mathbb R^n)$ satisfies
\begin{align}\label{ineq:mass-Gamma<=V}
\mathbb M(\Gamma)
\leq C\left(\mathbb M(V)+\mathcal E_p(V,\widetilde B)\right).
\end{align}
Moreover, if ${\rm spt}\|V\|\subset K$ for some compact $K\subset\mathbb R^n$, then, for every $p\in[1,\infty)$, there exists a positive constant $C_K=C_K(m,n,p,\gamma,K)$ such that
\begin{align}\label{ineq:mass-V<=Gamma}
\mathbb M(V)
\leq C_K\left(\mathbb M(\Gamma)+\mathcal E_p(V,\widetilde B)\right).
\end{align}
\end{lemma}

\begin{proof}
Since $\widetilde G_{m,n,\gamma}$ is compact and $\eta$ is continuous on $\widetilde G_{m,n,\gamma}$, by standard approximation there exists a compactly supported $C^1$-map $\widehat\eta:\mathbb R^n\times\Lambda_m(\mathbb R^n)\to\mathbb R^n$ such that
\begin{align*}
|\widehat\eta-\eta|
\leq\frac12\quad\text{on }\widetilde G_{m,n,\gamma},\quad
\|\widehat\eta\|_{C^1}
\leq C(m,n,\gamma).
\end{align*}

Let $\mathcal U_1(\gamma):=\left\{x\in\mathbb R^n:{\rm dist}(x,\gamma)<1\right\}$, let $\ell\in C_c^1(\mathbb R^n)$ be a standard cutoff function supported in $\mathcal U_1(\gamma)$, such that $0\leq\ell\leq1$, $\ell=1$ on a neighborhood of $\gamma$, and that $|D\ell|\leq C(\gamma)$.

Since $\eta$ is unit and $\ell=1$ on $\gamma$, we have
\begin{align*}
\int_{\widetilde G_{m,n,\gamma}}\left<\phi,\eta\right> d\Gamma
\geq\frac12\mathbb M(\Gamma).
\end{align*}
Testing \eqref{iden:ocv} with $\phi(x,\omega):=\ell(x)\widehat\eta(x,\omega)$ and using \eqref{eq:T-S-bound}, we get
\begin{align}
\mathbb M(\Gamma)
\leq C(m,n,\gamma)\left(V\left(\mathcal U_1(\gamma)\times\widetilde{\mathcal G}_{m,n}\right)+\int_{\mathcal U_1(\gamma)\times\widetilde{\mathcal G}_{m,n}}|\widetilde B(x,\omega)|\,dV(x,\omega)\right).
\end{align}
For $1\leq p<\infty$, since we have the pointwise inequality $t\leq1+t^p$ for any $t\geq0$, \eqref{ineq:mass-Gamma<=V} thus follows.

We next prove \eqref{ineq:mass-V<=Gamma}. Let $\chi\in C_c^1(\mathbb R^n)$ be such that $\chi=1$ on a neighborhood of $K\cup\gamma$, and test the curvature identity \eqref{iden:ocv} with $\phi(x,\omega):=\chi(x)x$.
Note that $D_x\phi={\rm Id}$ and $D_\omega\phi=0$ on ${\rm spt}\|V\|$, we thus obtain
\begin{align}
m\mathbb M(V)
=-\int_{\widetilde G_{m,n,\gamma}}\left<x,\eta(x,\omega)\right>d\Gamma(x,\omega)-\int\left<\mathcal T(\widetilde B)(x,\omega),x\right> dV(x,\omega),
\end{align}
where we have used \eqref{defn:trace-B} and $\left<{\rm Id},\mathbf q(\omega)\right>={\rm tr}\mathbf q(\omega)=m$. By \eqref{eq:T-S-bound}, we deduce
\begin{align}\label{ineq:mass-interme}
m\mathbb M(V)
\leq C_K\mathbb M(\Gamma)+C_K\int|\widetilde B|\,dV.
\end{align}
For $p=1$, this directly proves the assertion. For $p>1$, after increasing $C_K=C_K(m,n,p,\gamma,K)$, if necessary, we deduce by Young's inequality that
\begin{align*}
C_K|\widetilde B|
\leq\frac m2+\tilde C_K|\widetilde B|^p.
\end{align*}
Integrating against $V$, we get the mass term $\frac m2\mathbb M(V)$, which can be absorbed into the LHS of \eqref{ineq:mass-interme}, and \eqref{ineq:mass-V<=Gamma} is thus proved.
\end{proof}

\begin{proposition}\label{Prop:initial-cptness}
Let $p\in(1,\infty)$. For each $j\in\mathbb N$, let $(V_j,\widetilde B^j,\Gamma_j)$ satisfy Definition \ref{defn:ocv}, and suppose that $\sup_{j\in\mathbb N}\mathbb M(V_j)<\infty$. If
\begin{align}
\sup_{j\in\mathbb N}\mathcal E_p(V_j,\widetilde B^j)
\leq\Lambda<\infty.
\end{align}
Then, after passing to a subsequence, there exist $(V,\widetilde B,\Gamma)$ satisfying Definition \ref{defn:ocv}, such that $V_j\overset{*}{\rightharpoonup}V$ as oriented varifolds, $\Gamma_j\overset{*}{\rightharpoonup}\Gamma$ as Radon measures on $\widetilde G_{m,n,\gamma}$, and
\begin{align}\label{eq:lsc-L^p}
\mathcal E_p(V,\widetilde B)
\leq\liminf_{j\to\infty}\mathcal E_p(V_j,\widetilde B^j).
\end{align}
\end{proposition}

\begin{proof}
By Lemma \ref{Lem:bdry-mass-estimate} and the assumptions, we have $\sup_j\mathbb M(\Gamma_j)<\infty$. By compactness of Radon measures, after passing to a subsequence, $V_j\overset{*}{\rightharpoonup}V$ and $\Gamma_j\overset{*}{\rightharpoonup}\Gamma$ for some oriented varifold $V$ and Radon measure $\Gamma$ on $\widetilde G_{m,n,\gamma}$.

We regard $\widetilde B$ as a ${\rm Hom}(\mathbb R^n,\Lambda_m(\mathbb R^n))$-valued function, then define the vector-valued Radon measure $\mathcal B_j:=\widetilde B^jV_j$. By H\"older's inequality,
\begin{align*}
|\mathcal B_j|\left(\mathbb R^n\times\widetilde{\mathcal G}_{m,n}\right)
\leq\mathcal E_p(V_j,\widetilde B^j)^{1/p}\mathbb M(V_j)^{1-1/p}.
\end{align*}
Thus the total variations of $\mathcal B_j$ are uniformly bounded, and after passing to a further subsequence, $\mathcal B_j\overset{*}{\rightharpoonup}\mathcal B$ for some ${\rm Hom}(\mathbb R^n,\Lambda_m(\mathbb R^n))$-valued Radon measure $\mathcal B$.

By the standard compactness result for the so-called {\em measure-function pairs} (cf. \cite[Section 4]{Hutchinson86}), we see that $\mathcal B=BV$ for some $B\in L^p(V;({\rm Hom}(\mathbb R^n,\Lambda_m(\mathbb R^n)))$ and that \eqref{eq:lsc-L^p} holds. Moreover, for any fixed test function $\phi\in
C_c^1\left(\mathbb R^n\times\Lambda_m(\mathbb R^n);\mathbb R^n\right)$, the interior terms in \eqref{iden:ocv} for each $(V_j,\widetilde B_j,\Gamma_j)$ converge, while the boundary terms also converge thanks to the fact that $\eta$ is continuous. Letting $j\to\infty$, we thus obtain
\begin{align*}
\int_{\mathbb R^n\times\widetilde{\mathcal G}_{m,n}}\left(D_\omega\phi\cdot \widetilde{B}+\left<\mathcal T(\widetilde{B}),\phi\right>+\left< D_x\phi,\mathbf q(\omega)\right>\right) dV
=-\int_{\widetilde G_{m,n,\gamma}}\left<\phi,\eta\right> d\Gamma,
\end{align*}
which shows that $(V,\widetilde B,\Gamma)$ satisfies Definition
\ref{defn:ocv} and completes the proof.
\end{proof}

\subsection{Consequences of compactness}

\subsubsection{Minimization of the curvature energy}
We consider the following two minimization problems.
\begin{definition}\label{Defn:mmmzing}
\normalfont
For $p\in(1,\infty)$, define
\begin{align}\label{defn:kappa-inf}
\kappa_{m,\gamma}^p
:=\inf\left\{\mathcal E_p(V,\widetilde B):V\in\widetilde{\bf CV}_{m,\gamma}(\mathbb R^n), \mathbb M(V)=1\right\}.
\end{align}
\end{definition}

\begin{definition}\label{Defn:constrain-mmmzing}
\normalfont
For a compact domain $K\subset\mathbb R^n$ and $p\in(1,\infty)$, define
\begin{align}
\kappa_{m,\gamma,K}^p
:=\inf\left\{\mathcal E_p(V,\widetilde B):
V\in\widetilde{\bf CV}_{m,\gamma}(\mathbb R^n),
\mathbb M(V)=1,{\rm spt}\|V\|\subset K\right\}.
\end{align}
\end{definition}

\begin{remark}
\normalfont
The admissible classes in Definitions \ref{Defn:mmmzing} and \ref{Defn:constrain-mmmzing} are nonempty. Indeed, a round $m$-sphere contained in $\mathbb R^n$ (or in $K$) induces an oriented curvature varifold with boundary $\Gamma=0$, one can then normalize its mass to 1 and get a desired element in these admissible classes.
\end{remark}

\begin{proposition}\label{Prop:unconstrained-inf-0}
For every $p\in(1,\infty)$, one has
\begin{align*}
\kappa_{m,\gamma}^p=0.
\end{align*}
\end{proposition}

\begin{proof}
Let $\Sigma_r:=\mathbb S_r^m\subset\mathbb R^{m+1}\subset\mathbb R^n$ be a round sphere of radius $r$, and let $V_r$ be its associated oriented varifold normalized by $\mathbb M(V_r)=1$. The second fundamental form $A$ of $\Sigma_r$ satisfies $|A|=\frac{\sqrt m}{r}$. Also note that by \eqref{eq:widetilde-B=A-len} and Remark \ref{rem:scaling-prop}, the curvature field of $V_r$, say $\widetilde B_r$, satisfies $|\widetilde B_r|=|A|$. Hence for any $p\in(1,\infty)$,
\begin{align*}
\mathcal E_p(V_r,\widetilde B_r)
=m^{p/2}r^{-p}.
\end{align*}
Letting $r\to\infty$ proves the assertion.
\end{proof}

\begin{proposition}
Let $K\subset\mathbb R^n$ be a compact domain and $p\in(1,\infty)$, then the infimum $\kappa_{m,\gamma,K}^p$ is finite and is attained.
\end{proposition}

\begin{proof}
Let $\{(V_j,\widetilde B^j,\Gamma_j)\}_{j\in\mathbb N}$ be a minimizing sequence. By definition,
\begin{align*}
\mathbb M(V_j)=1,\quad{\rm spt}\|V_j\|\subset K,\quad\forall j\in\mathbb N.
\end{align*}
The curvature energies are uniformly bounded, and hence by
\ref{Lem:bdry-mass-estimate}, we get a uniform bound for
$\mathbb M(\Gamma_j)$. By Proposition \ref{Prop:initial-cptness}, we obtain a subsequential limiting triple $(V,\widetilde B,\Gamma)$ satisfying Definition \ref{defn:ocv}.
Since each $V_j$ is supported in $K$, we also have
\begin{align*}
\mathbb M(V)=1,\quad{\rm spt}\|V\|\subset K.
\end{align*}
By \eqref{eq:lsc-L^p} we see the infimum is attained by $(V,\widetilde B,\Gamma)$, which completes the proof.
\end{proof}

\subsubsection{Energy lower bound}

\begin{proposition}\label{Prop:constr-energy>0}
Suppose that the connected submanifold $\gamma$ is not contained in any affine $m$-plane. Then for any compact $K\subset\mathbb R^n$ with $\gamma\subset K$, and every $p\in(1,\infty)$,
\begin{align*}
\kappa_{m,\gamma,K}^p>0.
\end{align*}
\end{proposition}

\begin{proof}
We argue by contradiction and suppose that $(V_j,\widetilde B^j,\Gamma_j)$ satisfying Definition \ref{defn:ocv} is a sequence of triples such that $\mathbb M(V_j)=1$, ${\rm spt}\|V_j\|\subset K$, but
\begin{align*}
\mathcal E_p(V_j,\widetilde B^j)\to0\text{ as }j\to\infty.
\end{align*}
By Proposition \ref{Prop:initial-cptness}, we obtain a subsequential limiting triple $(V,\widetilde B,\Gamma)$ satisfying Definition \ref{defn:ocv} with $\mathbb M(V)=1$, ${\rm spt}\|V\|\subset K$, and $\widetilde B=0$ for $V$-a.e. However, this contradicts the structural theorem (Theorem \ref{Thm:structure-V-Gamma-S=0}) since $\gamma$ is not contained in an affine $m$-plane. The proof is thus completed.
\end{proof}

\section{Mass and topology bounds for immersions}\label{sec:6}

As shown in the smooth computations before Definition \ref{defn:ocv}, every $\Sigma\in\mathbf{S}^m(\gamma)$ induces a triple $(V,\widetilde B,\Gamma)$ satisfying Definition \ref{defn:ocv}, and by \eqref{eq:widetilde-B=A-len}, we have $|\widetilde B|=|A|$ where $A$ is the second fundamental form of the immersion.

\subsection{Mass bound in arbitrary dimension and codimension}\label{Sec:6-1}

\begin{proof}[Proof of Theorem \ref{Thm:area-bounds}]
Let $(V,\widetilde B,\Gamma)$ denote the triple induced by $\Sigma\in\mathbf{S}^m(\gamma)$. Put $\widetilde V:=\frac{1}{|\Sigma|}V$, $\widetilde\Gamma:=\frac{1}{|\Sigma|}\Gamma$, so that $\mathbb M(\widetilde V)=1$ and $(\widetilde V,\widetilde B,\widetilde \Gamma)$ satisfies Definition \ref{defn:ocv}, thanks to Remark \ref{rem:scaling-prop}.
By assumption there exists a compact set $K=K(\gamma,R)\subset\mathbb R^n$ such that ${\rm spt}\|\widetilde V_j\|\subset K$ for each $j\in\mathbb N$.
It follows that
\begin{align*}
\kappa_{m,\gamma,K}^p
\leq\mathcal E_p(\widetilde V,\widetilde B)
\leq\frac{1}{|\Sigma|}
\int_\Sigma|A|^p\,d\mathcal H^m,
\end{align*}
which shows the bound on $|\Sigma|$.
Taking also \eqref{ineq:mass-Gamma<=V} into account, we complete the proof.
\end{proof}

\subsection{Mass and topology bounds for surfaces}\label{sec:6-mass-topo-bds}

Throughout this subsection, let $\gamma\subset\mathbb R^n$ be a compact, connected, oriented, embedded $C^1$-curve, and let $\Sigma\in\mathbf{S}^2(\gamma)$ with immersion given by $F:\Sigma\to\mathbb R^n$ with the boundary immersion denoted as $\beta:=F|_{\partial\Sigma}$.
We write
\begin{align*}
\partial\Sigma
=S_1\cup\cdots\cup S_b
\end{align*}
as the disjoint union of the connected boundary components. Since $\gamma$ is connected and $\beta$ is an orientation-preserving local diffeomorphism, each restriction $\beta_a:=\beta|_{S_a}:S_a\to\gamma$ is a finite orientation-preserving covering. If we let
\begin{align*}
d_a:=\deg\beta_a\in\mathbb N
\end{align*}
and define the total boundary covering multiplicity by
\begin{align}
N(\beta):=\sum_{a=1}^b d_a,
\end{align}
then we have
\begin{align}\label{eq:bdry-length}
|\partial\Sigma|
=N(\beta)\mathcal H^1(\gamma),
\end{align}
and also $b\leq N(\beta)$.

For $N_0\in\mathbb N$, we consider the following subclass of $\mathbf{S}^2(\gamma)$:
\begin{align}\label{defn:surface-bdry-finite-cov}
\mathbf{S}^2_{N_0}(\gamma)
:=\left\{\Sigma\in\mathbf{S}^2(\gamma): \Sigma\text{ is connected and }N(\beta)\leq N_0\right\}.
\end{align}
The following estimate is a direct consequence of Simon's monotonicity formula, which we record in Appendix \ref{sec:Simon-monotonicity-bdry}.
\begin{lemma}\label{Lem:crtical-ener-gap}
Let $N_0\in\mathbb N$, and let $\Sigma\in\mathbf{S}^2_{N_0}(\gamma)$. If
\begin{align}\label{condi:8pi-threshold}
\int_\Sigma|A|^2\,d\mathcal H^2
<8\pi-\delta
\end{align}
for some $\delta\in(0,8\pi)$. Then
\begin{align}\label{ineq:dist-8pi-gap}
\sup_{x\in F(\Sigma)}{\rm dist}(x,\gamma)
\leq\frac{4N_0\mathcal H^1(\gamma)}{\delta},
\end{align}
and there exists $\mathcal R=\mathcal R(\gamma,\delta,N_0)<\infty$ such that $F(\Sigma)\subset\overline{B_{\mathcal R}(0)}$.
\end{lemma}

\begin{proof}
Let $x_0\in F(\Sigma)\setminus\gamma$. Then $x_0$ has an interior preimage under $F$, and $\Theta(F,x_0)\geq1$. By Corollary \ref{Cor:Simon-monotonicity-bdry}, one has
\begin{align}\label{ineq:pi<=dist}
\pi
\leq\frac1{16}\int_\Sigma|H|^2\,d\mathcal H^2+\frac{|\partial\Sigma|}{2{\rm dist}(x_0,\gamma).
}
\end{align}
Since $|H|^2\leq2|A|^2$, \eqref{condi:8pi-threshold} thus gives
\begin{align*}
\frac1{16}\int_\Sigma|H|^2\,d\mathcal H^2
<\pi-\frac{\delta}{8}.
\end{align*}
Using \eqref{eq:bdry-length} and $N(\beta)\leq N_0$ in
\eqref{ineq:pi<=dist}, we thus obtain
\begin{align*}
\frac{\delta}{8}
<\frac{N_0\mathcal H^1(\gamma)}{2{\rm dist}(x_0,\gamma)},
\end{align*}
which proves \eqref{ineq:dist-8pi-gap}.
The last assertion follows easily, which completes the proof.
\end{proof}

\begin{definition}\label{Defn:small-ener-class}
\normalfont
For $\delta\in(0,8\pi)$ and $N_0\in\mathbb N$, define the subclass of $\mathcal S^{N_0}_{2,\gamma}$ with small critical energy:
\begin{align}
\mathcal F_{\gamma,\delta,N_0}
:=\left\{\Sigma\in\mathbf{S}^2_{N_0}(\gamma):\int_\Sigma|A|^2\,d\mathcal H^2
<8\pi-\delta\right\}.
\end{align}
For $p\in(1,\infty)$, we consider the minimization problem
\begin{align}
\mathfrak h_{\gamma,\delta,N_0}^p
:=\inf\left\{\frac1{|\Sigma|}\int_\Sigma|A|^p\,d\mathcal H^2:\Sigma\in\mathcal F_{\gamma,\delta,N_0}\right\}.
\end{align}
\end{definition}

\begin{proposition}\label{Prop:area-energy-2d}
Suppose that $\gamma$ is not contained in any affine $2$-plane. Then, for every $p\in(1,\infty)$, $\delta\in(0,8\pi)$, and $N_0\in\mathbb N$, if $\mathcal F_{\gamma,\delta,N_0}\neq\emptyset$, then
\begin{align}
\mathfrak h_{\gamma,\delta,N_0}^p>0,
\end{align}
and every $\Sigma\in\mathcal F_{\gamma,\delta,N_0}$ satisfies
\begin{align}\label{esti:area-interior-bdry-L^p-energy-3}
|\Sigma|
\leq\frac1{\mathfrak h_{\gamma,\delta,N_0}^p}\int_\Sigma|A|^p\,d\mathcal H^2.
\end{align}
\end{proposition}

\begin{proof}
By Lemma \ref{Lem:crtical-ener-gap}, there exists $\mathcal R=\mathcal R(\gamma,\delta,N_0)$ such that for every $\Sigma\in\mathcal F_{\gamma,\delta,N_0}$, $F(\Sigma)$ is supported in $\overline{B_{\mathcal R}(0)}$. Let $(V,\widetilde B,\Gamma)$ denote the triple induced by $\Sigma$, and define the normalized oriented varifold and boundary measure by
\begin{align*}
\widetilde V
:=\frac1{|\Sigma|}V,\quad
\widetilde\Gamma
:=\frac1{|\Sigma|}\Gamma.
\end{align*}
Then $\mathbb M(\widetilde V)=1$ and $(\widetilde V,\widetilde B,\widetilde\Gamma)$ satisfies Definition \ref{defn:ocv} thanks to Remark \ref{rem:scaling-prop}, hence
\begin{align}
\kappa_{2,\gamma,\overline{B}_{\mathcal R}(0)}^p
\leq\mathcal E_p(\widetilde V,\widetilde B)
=\frac1{|\Sigma|}\int_\Sigma|A|^p \,d\mathcal H^2.
\end{align}
Taking infimum over $\mathcal F_{\gamma,\delta,N_0}$, we thus obtain
\begin{align*}
\mathfrak h_{\gamma,\delta,N_0}^p
\geq\kappa_{2,\gamma,\overline{B}_{\mathcal R}(0)}^p.
\end{align*}
Since $\gamma$ is not contained in an affine $2$-plane, the right-hand side is positive by Proposition \ref{Prop:constr-energy>0}. This proves the assertion.
\end{proof}

\begin{theorem}\label{Thm:area-geneus-surface}
Let $\gamma\subset\mathbb R^n$ be a compact, connected, oriented, embedded $C^2$ curve. Let $p\in[2,\infty)$, $N_0\in\mathbb N$, and let $\Sigma\in\mathbf{S}^2_{N_0}(\gamma)$. Assume that one of the following holds:
\begin{enumerate}
\item there exists a compact domain $K\subset\mathbb R^n$ such that $F(\Sigma)\subset K$;
\item there exists $\delta\in(0,8\pi)$ such that
\begin{align*}
\int_\Sigma|A|^2\,d\mathcal H^2
<8\pi-\delta.
\end{align*}
\end{enumerate}
Then we have:
\begin{enumerate}
\item [(i)] If $\gamma$ is not contained in an affine $2$-plane, then
\begin{align}\label{esti:Euler-Lp-energy}
g(\Sigma)
\leq C\left(1+\int_\Sigma|A|^p \,d\mathcal H^2\right),
\end{align}
where $C$ depends on $p,\gamma,N_0$ and on $K$ in case $(1)$, or on $\delta$ in case $(2)$.

\item [(ii)] If $\gamma$ is contained in an affine $2$-plane, then
\begin{align}\label{ineq:planar-area-topology}
|\Sigma|
\leq C\left[\int_\Sigma|A|^2\,d\mathcal H^2+2\pi\left|2-2g(\Sigma)-b(\Sigma)\right|\right],
\end{align}
where $C$ depends on $\gamma,N_0$ and on $K$ in case $(1)$, or on
$\delta$ in case $(2)$.
\end{enumerate}
\end{theorem}

\begin{proof}
Let $b$ denote the number of boundary connected components of $\Sigma$, then by definition $b\leq N(\beta)\leq N_0$, and since $\Sigma$ is connected and oriented, one has
\begin{align}\label{eq:Euler-Chara}
\chi(\Sigma)
=2-2g(\Sigma)-b.
\end{align}
By the Gauss-Bonnet theorem,
\begin{align}\label{eq:Gauss-Bonnet}
2\pi\chi(\Sigma)
=\int_\Sigma K\,d\mathcal H^2+\int_{\partial\Sigma}\kappa_{g}\,d\mathcal H^1,
\end{align}
where $K$ denotes the Gauss curvature of $\Sigma$ and $\kappa_{g}(z):=\left< (D_{\tau_\gamma}\tau_\gamma)(F(z)),\eta_{\partial\Sigma}(z)\right>$ denotes the geodesic curvature of $\partial\Sigma$. We note that this is an immediate consequence if $\Sigma$ is of class $C^2$ up to the boundary, and for general $\Sigma\in\mathbf{S}^2(\gamma)$ this follows by approximation, see \cite[Proof of Theorem 1]{KM22}.
We also note that
\begin{align}\label{eq:Gauss-2nd-fdmtal-form}
|K|
\leq\frac12|A|^2,
\end{align}
since $\Sigma$ is $2$-dimensional.
Moreover, along $\partial\Sigma$,
\begin{align*}
|\kappa_{g}|
\leq|D_{\tau_\gamma}\tau_\gamma|\circ F.
\end{align*}
Since $\beta$ is a covering of $\gamma$, we get
\begin{align}
\int_{\partial\Sigma}|\kappa_{g}| \,d\mathcal H^1
&\leq N(\beta)\int_\gamma|D_{\tau_\gamma}\tau_\gamma|\,d\mathcal H^1
\leq N_0\int_\gamma|D_{\tau_\gamma}\tau_\gamma|\,d\mathcal H^1.
\end{align}
Back to \eqref{eq:Gauss-Bonnet} we thus find
\begin{align*}
2\pi|\chi(\Sigma)|
\leq\frac12\int_\sigma|A|^2\,d\mathcal H^2+N_0\int_\gamma|D_{\tau_\gamma}\tau_\gamma|\,d\mathcal H^1.
\end{align*}
Using \eqref{eq:Euler-Chara} we see $g(\Sigma)=(2-\chi(\Sigma)-b)/2\leq1+|\chi(\Sigma)|/2$, it follows that
\begin{align}\label{eq:genus-L2-0}
g(\Sigma)
\leq C(\gamma,N_0)\left(1+\int_\Sigma|A|^2\,d\mathcal H^2\right),
\end{align}
which proves \eqref{esti:Euler-Lp-energy} for the case $p=2$. For $p>2$, H\"older's inequality gives
\begin{align}\label{eq:L2-from-Lp-area}
\int_\Sigma|A|^2\,d\mathcal H^2
\leq\left(\int_\Sigma|A|^p \,d\mathcal H^2\right)^{2/p}|\Sigma|^{1-2/p}.
\end{align}
Under condition $(1)$, we can apply Theorem \ref{Thm:area-bounds} to get
\begin{align*}
|\Sigma|
\leq C\int_\Sigma|A|^p \,d\mathcal H^2.
\end{align*}
Under condition $(2)$, the same estimate follows from Proposition
\ref{Prop:area-energy-2d}. In both cases, back to \eqref{eq:L2-from-Lp-area} we thus obtain
\begin{align*}
\int_\Sigma|A|^2\,d\mathcal H^2
\leq C\int_\Sigma|A|^p \,d\mathcal H^2.
\end{align*}
Combining with \eqref{eq:genus-L2-0}, \eqref{esti:Euler-Lp-energy} is thus proved.

Now we consider the case that $\gamma$ is contained in an affine $2$-plane. To prove \eqref{ineq:planar-area-topology}, we argue by contradiction and suppose that \eqref{ineq:planar-area-topology} fails. Then there exists a sequence $\Sigma_j\in\mathbf{S}^2_{N_0}(\gamma)$ (with the immersion map $F_j:\Sigma_j\to\mathbb R^n$) satisfying the same constraint as in the statement, but with
\begin{align}\label{eq:planar-contra-assum}
\frac{\displaystyle\int_{\Sigma_j}|A|^2\,d\mathcal H^2+2\pi\left|2-2g(\Sigma_j)-b(\Sigma_j)\right|}{|\Sigma_j|}\to0.
\end{align}

Under condition $(1)$, all $F_j(\Sigma_j)$ lie in $K$. Under condition $(2)$, we can again find such a compact set $K$ thanks to Lemma \ref{Lem:crtical-ener-gap}. Thus in both cases, there exists a compact set $K\subset\mathbb R^n$, such that
\begin{align*}
F_j(\Sigma_j)\subset K \quad\forall j\in\mathbb N.
\end{align*}

Let $(V_j,\widetilde B_j,\Gamma_j)$ denote the triple induced by $\Sigma_j$, and define the normalized oriented varifolds and boundary measures by
\begin{align*}
\widetilde V_j
:=\frac1{|\Sigma_j|}V_{j},\quad
\widetilde\Gamma_j
:=\frac1{|\Sigma_j|}\Gamma_{j}.
\end{align*}
Then $\mathbb M(\widetilde V_j)=1$, and
\begin{align}
\mathcal E_2(\widetilde V_j,\widetilde B_{j})
\leq\frac1{|\Sigma_j|}\int_{\Sigma_j}|A|^2\,d\mathcal H^2
\to0.
\end{align}
By Proposition \ref{Prop:initial-cptness}, after passing to a subsequence, there exist $V$ and $\Gamma$ such that $\widetilde V_j\to V$ as oriented varifolds, $\widetilde\Gamma_j\overset{*}{\rightharpoonup}\Gamma$, and $(V,0,\Gamma)$ satisfies Definition \ref{defn:ocv}.
Moreover, we have $\mathbb M(V)=1$ and ${\rm spt}\|V\|\subset K$.

To proceed, we define on $\widetilde G_{1,2,\gamma}$ the function
\begin{align}
\mathfrak k_\gamma(x,\omega)
:=\left< D_{\tau_\gamma}\tau_\gamma(x),\eta(x,\omega)\right>,
\end{align}
then consider a compactly supported continuous extension $\widehat{\mathfrak k}_\gamma:\mathbb R^n\times\Lambda_2(\mathbb R^n)\to\mathbb R$ of $\mathfrak k_\gamma$. Since $\partial\Sigma$ is an oriented covering of $\gamma$, we have
\begin{align}
\widetilde\Gamma_j\left(\widehat{\mathfrak k}_\gamma\right)
=\frac1{|\Sigma_j|}\int_{\partial\Sigma_j}\kappa_{g}\,d\mathcal H^2.
\end{align}
By \eqref{eq:Gauss-Bonnet}, \eqref{eq:Gauss-2nd-fdmtal-form} and \eqref{eq:planar-contra-assum}, we find
\begin{align}
\left|\widetilde\Gamma_j\left(\widehat{\mathfrak k}_\gamma\right)\right|
\leq\frac{2\pi\left|2-2g(\Sigma_j)-b(\Sigma_j)\right|}{|\Sigma_j|}+\frac1{2|\Sigma_j|}\int_{\Sigma_j}|A|^2\,d\mathcal H^2
\to0,
\end{align}
which implies that
\begin{align}\label{eq:limit-Gamma-vani}
\Gamma\left(\widehat{\mathfrak k}_\gamma\right)
=0.
\end{align}

Let $L$ be the affine $2$-plane containing $\gamma$, and let $\Delta$ be the unique bounded connected component of $L\setminus\gamma$. Using Theorem
\ref{Thm:structure-V-Gamma-S=0} for the zero-curvature varifold $V$, if we put
\begin{align}
W
:=\int_{\widetilde{\mathcal G}_{2,n}}\int_{\mathbb P_\omega^\perp}\sum_{\Delta'\in\mathscr S_{2,\gamma}(L_{\omega,z})}\theta_{\omega,z,\Delta'}d\lambda_\omega(z)d\nu(\omega),
\end{align}
then by \eqref{eq:V-struc} we have
\begin{align*}
1
=\mathbb M(V)
=W\mathcal H^2(\Delta),
\end{align*}
so that $W=\frac1{\mathcal H^2(\Delta)}>0$.

On the other hand, using Corollary \ref{Cor:structure-bdry} with the test function $\mathfrak k_\gamma$, we obtain
\begin{align}
\Gamma\left(\widehat{\mathfrak k}_\gamma\right)
&=\int_{\widetilde{\mathcal G}_{2,n}}\int_{\mathbb P_\omega^\perp}\sum_{\Delta'\in\mathscr S_{2,\gamma}(L_{\omega,z})}\theta_{\omega,z,\Delta'}\\
&\quad\cdot\int_\gamma\left<\eta_\Delta^{\mathrm{in}}(x),\eta(x,\omega) \right>\left<D_{\tau_\gamma}\tau_\gamma(x),\eta(x,\omega)\right> d\mathcal H^1(x)d\lambda_\omega(z)d\nu(\omega).
\end{align}
Note that both $\eta_\Delta^{\mathrm{in}}(x)$ and $\eta(x,\omega)$ are unit vectors in the one-dimensional space $L\cap \left(T_x\gamma^\perp\right)$,
therefore $\eta(x,\omega)=\pm\eta_\Delta^{\mathrm{in}}(x)$, and hence
\begin{align}
\left<\eta_\Delta^{\mathrm{in}},\eta(x,\omega) \right>\left<D_{\tau_\gamma}\tau_\gamma,\eta(x,\omega) \right>
=\left< D_{\tau_\gamma}\tau_\gamma,\eta_\Delta^{\mathrm{in}} \right>.
\end{align}
By the Hopf Umlaufsatz,
\begin{align*}
\int_\gamma\left<D_{\tau_\gamma}\tau_\gamma,\eta_\Delta^{\mathrm{in}}\right>\,d\mathcal H^1
=2\pi,
\end{align*}
and it follows that
\begin{align}
\Gamma\left(\widehat{\mathfrak k}_\gamma\right)
=2\pi W
=\frac{2\pi}{\mathcal H^2(\Delta)}>0,
\end{align}
which contradicts \eqref{eq:limit-Gamma-vani} and completes the proof.
\end{proof}
}

\section{Existence of minimizers of the boundary value problem with obstacle}\label{sec:7}

\subsection{Integral oriented varifolds}
In \cite{Hutchinson86}, the notion of so-called integral oriented varifolds is introduced, which we briefly recall.
\begin{definition}\label{Defn:oriented-integral-vflds}
\normalfont
Given an \textit{oriented $m$-varifold} $V$ on $U\subset\mathbb R^n$.
$V$ is called an \textit{oriented integral $m$-varifold} on $U$, if there exist an $m$-rectifiable set $\mathcal R$ in $U$, a measurable choice of orientation $\xi:\mathcal R\to\widetilde G(m,n)$ such that for $\mathcal H^m$-a.e. $x\in\mathcal R$, $\xi(x)\in\widetilde G(m,n)$ is one of the two choices of unit simple $m$-vector associated to the approximate tangent space $T_x\mathcal R$, and a couple $(\theta_1,\theta_2)$ of locally integrable $\mathbb N$-valued functions on $\mathcal R$, with $(\theta_1,\theta_2)\neq(0,0)$ for $\mathcal H^m$-a.e. $x\in\mathcal R$, so that $V$ acts on $\phi\in C_c^0(U\times\Lambda_m(\mathbb R^n))$ as follows:
\begin{align}\label{defn:V-OIV}
V(\phi)
=\int_\mathcal R\theta_1(x)\phi(x,\xi(x))+\theta_2(x)\phi(x,-\xi(x))\,{d}\mathcal H^m(x).
\end{align}
Namely, $V=\mathcal H^m\llcorner\mathcal R\otimes\left(\theta_1(x)\delta_{\xi(x)}+\theta_2(x)\delta_{-\xi(x)}\right)$.
For simplicity we write $V={\mathbf v}(\mathcal R,\xi,\theta_1,\theta_2)$.
\end{definition}

For any oriented $m$-varifold $V$ in $U$, the $m$-current $\mathbf{c}(V)$ associated to it is defined by its action on an arbitrary $m$-form $\alpha$ having compact support in $U$ as follows:
\begin{align}\label{defn:c(V)(alpha)}
\mathbf{c}(V)(\alpha)
=\int\left<\omega,\alpha(x)\right>{d} V(x,\omega).
\end{align}
In particular, for $V=\mathbf{v}(\mathcal{R},\xi,\theta_1,\theta_2)$ we can write explicitly
\begin{align}\label{eq:indeuced-current-ori-rec-vfld}
\mathbf{c}(V)(\alpha)
=\int_\mathcal{R}(\theta_1(x)-\theta_2(x))\left<\xi(x),\alpha(x)\right>{d}\mathcal H^{m}(x),
\end{align}
which agrees with Hutchinson's notion of induced current \cite[pp. 50]{Hutchinson86}.

We are going to use the following celebrated compactness theorem for oriented integral varifolds by Hutchinson \cite[Theorem 3.1]{Hutchinson86}.
\begin{theorem}\label{Thm:Hutchinson-cpt}
Let $K\subset\mathbb R^n$ be a compact set.
Let $\{V_j\}_{j\in\mathbb N}$ be oriented integral $m$-varifolds on $K$ such that
\begin{align}\label{eq:Hutchison-bound}
\sup_{j\in\mathbb N}\left\{\|{\mathbf q_\#V_j}\|(K)+\|{\delta(\mathbf q_\#V_j)}\|(K)+\|(\partial\mathbf{c}(V_j))\|(K)\right\}<\infty.
\end{align}
Then there exists an oriented integral $m$-varifold $V$ on $K$, to which $V_j$ subsequentially converges, as oriented varifolds.
\end{theorem}

\subsection{Minimizers of the boundary value problem with obstacle}


To minimize the $L^p$-curvature energy under the given boundary condition, we are in particular interested in the integer rectifiable subclass of $\widetilde{\bf CV}_{m,\gamma}(\mathbb R^n)$, say $\widetilde{\bf ICV}_{m,\gamma}(\mathbb R^n)$, whose elements are such that $V\in\widetilde{\bf CV}_{m,\gamma}(\mathbb R^n)$ with $V$ an integral oriented $m$-varifold. 

We first use the curvature identity \eqref{iden:ocv} to estimate the boundary mass of the induced current $\mathbf c(V)$.
\begin{lemma}
Let $K\subset\mathbb R^n$ be a compact set. Let $(V,\widetilde B,\Gamma)$ satisfy Definition \ref{defn:ocv} with $V$ integral, and assume that ${\rm spt}\|V\|\subset K$, then
\begin{align}\label{esti:bdry-current}
\|(\partial\mathbf{c}(V))\|(K)
\leq C(m,n)\int|\widetilde B|\,dV+\mathbb M(\Gamma).
\end{align}
\end{lemma}
\begin{proof}

For any smooth compactly supported $(m-1)$-form $\alpha$ on
$\mathbb R^n$, by definition we have
\begin{align}\label{eq:partial-c(V)}
\partial\mathbf{c}(V)(\alpha)
=\mathbf{c}(V)(d\alpha)
\overset{\eqref{defn:c(V)(alpha)}}{=}\int\left<\omega,d\alpha(x)\right>dV(x,\omega).
\end{align}

Now we choose any smooth cutoff $\chi\in C^\infty_c(\Lambda_m(\mathbb R^n))$ which is identically $1$ near $\widetilde{\mathcal G}_{m,n}$, and define a compactly supported test function $\Phi_\alpha$ by letting
\begin{align}\label{defn:Phi_alpha}
\left<\Phi_\alpha(x,\omega),v\right>
:=\chi(\omega)\left<\iota_v\omega,\alpha(x)\right>,\quad\forall v\in\mathbb R^n.
\end{align}
Fix any $(x,\omega)\in\mathbb R^n\times\widetilde{\mathcal G}_{m,n}$, we can write locally $\omega=e_1\wedge\cdots\wedge e_m$, where $e_1,\cdots,e_m$ is orthonormal, and compute
\begin{align}
\left<D_x\Phi_\alpha,\mathbf q(\omega)\right>
=&\left<D_x\Phi_\alpha,\sum_{i=1}^me_i\otimes e_i\right>
=\sum_{i=1}^m\left<D_{e_i}\Phi_\alpha,e_i\right>
=\sum_{i=1}^m\left<\iota_{e_i}\omega,D_{e_i}\alpha\right>\\
=&\sum_{i=1}^m(-1)^{i-1}(D_{e_i}\alpha)(e_1,\cdots,\widehat{e}_i,\cdots,e_m)
=d\alpha(e_1,\cdots,e_m)
=\left<\omega,d\alpha\right>,
\end{align}
where we have used the fact that $\iota_{e_i}\omega=(-1)^{i-1}e_1\wedge\cdots\wedge\widehat{e}_i\wedge\cdots\wedge e_m$.
Putting this back into \eqref{eq:partial-c(V)} and then using the curvature identity \eqref{iden:ocv}, we thus obtain
\begin{align}
\partial\mathbf{c}(V)(\alpha)
=&\int\left<D_x\Phi_\alpha,\mathbf q(\omega)\right>dV\\
=&-\int\left(D_\omega\Phi_\alpha\cdot\widetilde{B}+\left<\mathcal T(\widetilde{B}),\Phi_\alpha\right>\right)dV-\int\left<\Phi_\alpha,\eta\right> d\Gamma.
\end{align}

Now for every $(x,\omega)\in\mathbb R^n\times\widetilde{\mathcal G}_{m,n}$ and every $\xi\in\Lambda_m(\mathbb R^n)$, by taking derivative with respect to $t$ of \eqref{defn:Phi_alpha} along $\omega+t\xi$, we find
\begin{align}
    D_\omega\Phi_\alpha(x,\omega)[\xi]
    =\sum_{i=1}^n\left<\iota_{e_i}\xi,\alpha(x)\right>e_i,
\end{align}
where we have extended $\{e_1,\cdots,e_m\}$ to an orthonormal basis $\{e_1,\cdots,e_n\}$ of $\mathbb R^n$.
Thus it is easy to see that
\begin{align*}
|\Phi_\alpha|+|D_\omega\Phi_\alpha|
\leq C(m,n)|\alpha(x)|,\quad\forall(x,\omega)\in\mathbb R^n\times\widetilde{\mathcal G}_{m,n}.
\end{align*}
Taking $|\mathcal T(\widetilde B)|\leq C(m,n)|\widetilde B|$ into account, we deduce
\begin{align*}
|\partial\mathbf{c}(V)(\alpha)|
\leq C(m,n)\|\alpha\|_{C^0}\left(\int|\widetilde B|\,dV+\mathbb M(\Gamma)\right),
\end{align*}
which then yields the desired estimate \eqref{esti:bdry-current}.
\end{proof}

\begin{proof}[Proof of Theorem \ref{thm:minimizer-exis}]
Let $(V_j,\widetilde B_j,\Gamma_j)$ satisfy Definition \ref{defn:ocv} with $V_j$ integral be a minimizing sequence of \eqref{defn:mmmzer-L^p}. 

Since $\mathbb M(\Gamma_j)=k\mathcal H^{m-1}(\gamma)$, by \eqref{ineq:mass-V<=Gamma} we find
\begin{align*}
\sup_{j\in\mathbb N}\left\{\|V_j\|(\overline{B_R(0)})\right\}<\infty.
\end{align*}
By H\"older inequality, for the minimizing sequence, there holds
\begin{align*}
\sup_{j\in\mathbb N}\left\{\int|\widetilde B_j|\,dV_j\right\}
\leq\sup_{j\in\mathbb N}\left\{\left(\int|\widetilde B_j|^p\,dV_j\right)^{1/p}\|V_j\|(\overline{B_R(0)})^{1-1/p}\right\}
<\infty.
\end{align*}
By \eqref{esti:1st-vari} and \eqref{esti:bdry-current}, we also have
\begin{align*}
&\sup_{j\in\mathbb N}\left\{\|\delta(\mathbf q_\#V_j)\|(\overline{B_R(0)})
+\|\partial\mathbf c(V_j)\|(\overline{B_R(0)})\right\}\\
\leq&C(m,n)\sup_j\left\{\int|\widetilde B_j|\,dV_j\right\}+2k\mathcal H^{m-1}(\gamma)
<\infty.
\end{align*}
Thus we can apply the integral compactness Theorem \ref{Thm:Hutchinson-cpt}, together with Proposition \ref{Prop:initial-cptness}, to obtain a subsequential limit (still indexed by $j\in\mathbb N$) $(V,\widetilde B,\Gamma)$ satisfying Definition \ref{defn:ocv} with $V$ integral, such that $V_j\to V$ as varifolds, $\Gamma_j\to\Gamma$ as Radon measures, with ${\rm spt}\|V\|\subset\overline{B_R(0)}$, $\pi_\#\Gamma=k\mathcal H^{m-1}\llcorner\gamma$, as well as
\begin{align*}
\int|\widetilde B|^p\,dV
\leq\liminf_{j\to\infty}\int|\widetilde B_j|^p\,dV_j
=\beta^p_{k,R}(\gamma).
\end{align*}
Since the limiting varifold $V$ belongs to the same class in \eqref{defn:mmmzer-L^p}, we thus deduce that $\beta^p_{k,R}(\gamma)$ is attained.

Finally, if $\gamma$ is not contained in any affine
$m$-plane, then a triple $(V,\widetilde B,\Gamma)$ that attains $\beta^p_{k,R}(\gamma)=0$ would satisfy $\widetilde B=0$ for $V$-a.e.
By Theorem \ref{Thm:structure-V-Gamma-S=0} we must have $V=0$, and hence by Proposition \ref{Prop:bdry-spt}, it must be that $\pi_\#\Gamma=0$. However, this contradicts
$\pi_\#\Gamma=k\mathcal H^{m-1}\llcorner\gamma\ne0$ and hence $\beta^p_{k,R}(\gamma)>0$.
This completes the proof.
\end{proof}

\appendix

{

\section{Simon's monotonicity formula for immersed surfaces with boundary}\label{sec:Simon-monotonicity-bdry}

Let $F:\Sigma\longrightarrow\mathbb R^n$ be a compact, oriented, $C^1$ immersion of a surface with nonempty boundary, which is $C^2$ in the interior. Let $H$ denote the mean curvature vector and let $\eta$ denote the inward-pointing unit conormal along $\partial\Sigma$. By the tangential divergence theorem, for any $X\in C^1(\Sigma;\mathbb R^n)$,
\begin{align}\label{eq:1st-variation-app}
\int_\Sigma{\rm div} X\,d\mathcal H^2
=-\int_\Sigma\left< H,X\right>
\,d\mathcal H^2-\int_{\partial\Sigma}
\left<\eta,X\right> \,d\mathcal H^1.
\end{align}

Fix $x_0\in\mathbb R^n\setminus F(\partial\Sigma)$, for simplicity we adopt the following notations. For $y\in\Sigma$, set $r(y):=|F(y)-x_0|$, and for $\rho>0$, define
\begin{align*}
\Sigma_\rho(x_0)
:=\left\{y\in\Sigma:r(y)<\rho\right\}.
\end{align*}
We also define the density function as
\begin{align*}
\Theta(F,x_0)
:=\#F^{-1}(\{x_0\})
=\lim_{\rho\searrow0}\frac{|\Sigma_\rho(x_0)|}{\pi\rho^2}.
\end{align*}

The following Simon-type monotonicity formula is well-known to experts, cf. \cite{Simon93,KS04}, and \cite{NP20} for the version involving boundary terms.

\begin{proposition}\label{Prop:Simon-monotonicity-bdry}
Fix $x_0\in\mathbb R^n\setminus F(\partial\Sigma)$, for $\rho>0$, define
\begin{align}
\mathfrak g_{x_0}(\rho)
:=&\frac{|\Sigma_\rho(x_0)|}{\rho^2}+\frac1{16}\int_{\Sigma_\rho(x_0)}|H|^2\,d\mathcal H^2\\
&+\frac1{2\rho^2}\int_{\Sigma_\rho(x_0)}\left< H,F(y)-x_0 \right> \,d\mathcal H^2(y)\\
&-\frac12\int_{\partial\Sigma\cap\Sigma_\rho(x_0)}\left(\frac1{|F(y)-x_0|^2}-\frac1{\rho^2}\right)\left< F(y)-x_0,\eta \right>d\mathcal H^1(y).
\end{align}
Then $\mathfrak g_{x_0}$ is nondecreasing. In fact, for every $0<\rho_1<\rho_2$,
\begin{align}\label{eq:Simon-monotone}
\mathfrak g_{x_0}(\rho_1)+\int_{\Sigma_{\rho_2}(x_0)\setminus\Sigma_{\rho_1}(x_0)}\left|\frac{H}{4}+\frac{(F(y)-x_0)^\perp}{|F(y)-x_0|^2}\right|^2d\mathcal H^2(y)
=\mathfrak g_{x_0}(\rho_2).
\end{align}
Here $(F(\cdot)-x_0)^\perp$ denotes the orthogonal projection of $F(\cdot)-x_0$ onto the normal space of the immersion.
\end{proposition}

\begin{proof}[Sketch of proof]
Testing \eqref{eq:1st-variation-app} with the radial vector fields of the form $X(y)=\zeta(r(y))\left(F(y)-x_0\right)$, where $\zeta$ is a cut-off function, yields the results.
\end{proof}

\begin{corollary}\label{Cor:Simon-monotonicity-bdry}
Let $x_0\in F(\Sigma\setminus\partial\Sigma)\setminus F(\partial\Sigma)$.
Then
\begin{align}\label{eq:global-Simon-mono-bdry}
&\pi\Theta(F,x_0)+\int_\Sigma\left|\frac{H}{4}+\frac{(F-x_0)^\perp}{|F-x_0|^2}\right|^2d\mathcal H^2\\
=&\frac1{16}\int_\Sigma|H|^2\,d\mathcal H^2-\frac12\int_{\partial\Sigma}\frac{\left< F(y)-x_0,\eta \right>}{|F(y)-x_0|^2}\,d\mathcal H^1(y).
\end{align}
Consequently,
\begin{align}\label{ineq:Simon-monotone}
\pi\Theta(F,x_0)
\leq\frac1{16}\int_\Sigma|H|^2\,d\mathcal H^2+\frac{|\partial\Sigma|}{2{\rm dist}\left(x_0,F(\partial\Sigma)\right)}.
\end{align}
\end{corollary}

\begin{proof}
Let $\rho_1\searrow0$ and $\rho_2\to\infty$ in \eqref{eq:Simon-monotone}. Since $x_0$ is not contained in the boundary image, the boundary term vanishes for all sufficiently small $\rho_1$. Since the immersion is $C^2$ away from the boundary, we clearly have
\begin{align*}
\lim_{\rho\to0}\frac{|\Sigma_\rho(x_0)|}{\rho^2}
=\pi\Theta(F,x_0),
\end{align*}
and the remaining terms vanish as $\rho\to0$.

Since $\Sigma$ is compact, as $\rho\to\infty$ the area-ratio term and the term containing $\rho^{-2}$ in \eqref{eq:Simon-monotone} vanish, while the boundary term converges to
\begin{align*}
-\frac12\int_{\partial\Sigma}\frac{\left<F(y)-x_0,\eta \right>}{|F(y)-x_0|^2}\,d\mathcal H^1(y).
\end{align*}
This proves \eqref{eq:global-Simon-mono-bdry}.
Note that $\left|\frac{\left< F(y)-x_0,\eta \right>}{|F(y)-x_0|^2}\right|\leq\frac1{{\rm dist}\left(x_0,F(\partial\Sigma)\right)}$, \eqref{ineq:Simon-monotone} then follows.
\end{proof}
}


\bibliography{BibTemplate.bib}
\bibliographystyle{amsalpha}

\end{document}